\documentclass[12pt]{article}
\usepackage{amsmath,amssymb}
\usepackage{amsthm}
\usepackage{calc}
\usepackage{epsfig,afterpage}
\advance\textwidth by +1.0in
\advance\textheight by +1.0in
\advance\oddsidemargin by -0.5in
\advance\evensidemargin by -1.0in
\advance\topmargin by -0.5in
\def\be{\begin{equation}}
\def\ee{\end{equation}}
\def\bal{\begin{aligned}}
\def\eal{\end{aligned}}
\def\bq{\begin{eqnarray}}
\def\eq{\end{eqnarray}}
\def\beq{\begin{eqnarray*}}
\def\eeq{\end{eqnarray*}}
\def\ve{{\varepsilon}}
\def\btab{\begin{tabular}}
\def\etab{\end{tabular}}
\def\ba{\begin{array}}
\def\ea{\end{array}}
\def\bth{\begin{theorem}}
\def\eth{\end{theorem}}
\def\blm{\begin{lemma}}
\def\elm{\end{lemma}}
\def\bdf{\begin{definition}}
\def\edf{\end{definition}}
\def\bpr{\begin{proposition}}
\def\epr{\end{proposition}}
\def\brm{\begin{remark}}
\def\erm{\end{remark}}
\def\bnot{\begin{notation}}
\def\enot{\end{notation}}
\def\bobs{\begin{observation}}
\def\eobs{\end{observation}}
\def\bcrl{\begin{corrolary}}
\def\ecrl{\end{corrolary}}

\newcommand{\ab}{\mbox{\boldmath $a$}}

\newcommand{\const}{{\mbox{constant}}}
\newcommand{\sign}{\mbox{\rm sign\,}}
\newcommand{\Res}{\mbox{\rm Res\,}}
\newcommand{\Discrim}{\mbox{\rm Discrim\,}}
\newcommand{\Jacob}{\mbox{\rm Jacob\,}}
\newcommand{\Hess}{\mbox{\rm Hess\,}}
\newcommand{\Aff}{\mbox{\it Af\mbox{}f\,}}
\newcommand{\Supp}{\mbox{\rm Supp\,}}

\newcommand{\SSS}{{\bf S }}

\newcommand{\QS}{\mbox{\bf QS }}

\newcommand{\PP}{{\mathbb P}}
\newcommand{\R}{\mathbb{R}}
\newcommand{\C}{\mathbb{C}}

\newcommand{\N}{\mathbb{N}}
\newcommand{\Z}{\mathbb{Z}}
\newcommand{\AB}{{\bf A}}
\newcommand{\LL}{{\bf L}}
\newcommand{\EProof}{\ \hfill\rule[-1mm]{2mm}{2.5mm}}
\newcommand{\BProof}{\noindent{\it Proof:}\ \,}

\newtheorem{theorem}{Theorem}[section]
\newtheorem{lemma}{Lemma}[section]
\newtheorem{definition}{Definition}[section]
\newtheorem{proposition}{Proposition}[section]
\newtheorem{remark}{Remark}[section]
\newtheorem{notation}{Notation}[section]
\newtheorem{observation}{Observation}[section]
\newtheorem{corrolary}{Corrolary}[section]
\numberwithin{equation}{section}

\begin{document}
\setcounter{secnumdepth}{5}

\title{ Geometry of quadratic differential systems in the
neighbourhood of   infinity}

\author{Dana SCHLOMIUK\thanks{Work supported by NSERC and
        by the Quebec Education Ministry}\\
        D\'epartement de Math\'ematiques et de Statistiques\\
        Universit\'e de Montr\'eal
 \and
        Nicolae VULPE\thanks{Partially supported by NSERC}\\
        Istitute of Mathematics and Computer Science\\
        Academy of Science of Moldova }
\date{}

\maketitle
\begin{abstract}
In this article we consider the behavior in the vicinity of
infinity of the class of all planar quadratic differential
systems. This family depends on twelve parameters but due to
action of the affine group and re-scaling of time the family
actually depends on five parameters. We give simple,
integer-valued geometric invariants for this group action which
classify this family according to the topology of their phase
portraits in the vicinity of infinity. For each one of the classes
obtained we give necessary and sufficient conditions  in terms of
algebraic invariants and comitants so as to be able to easily
retrieve for any system, in any chart, the geometric as well as
the dynamic characteristics of the systems in the neighborhood of
infinity. A program was implemented for computer calculations.
\medskip

\noindent {\bf Keywords:} Poincar\'e compactification, singular
point, phase portrait, topological index,  intersection
multiplicity, linear group, affine invariant.
\end{abstract}
\section{ Introduction}
 We consider real planar polynomial differential system,
i.e.  systems  of  the  form
$$
\frac{dx}{dt} = p(x,y),\quad \frac{dy}{dt} = q(x,y)\eqno(S)
$$
where  $p$  and  $q$  are  polynomials  in  $x$  and  $y$  with real
coefficients  ($p,q\in \mathbb R[x,y]$).  In  this  article, a system of the above
form  with  $\max(\deg(p),\deg(q)) = 2$ will be called quadratic.

These are the simplest nonlinear  differential  systems.  However,  global
problems regarding  this  class are difficult to solve. In 1900  Hilbert  gave
his list of 23 problems and one of the few of them  still unsolved, the second
part of Hilbert's 16th problem, is  on planar polynomial systems. This problem
which asks for  the  maximum  H(n)  of the numbers  of limit cycles occurring
in  differential systems with $\max(\deg(p),\deg(q))\le n$ (for a discussion of
this problem cf. \cite{Dana_Sab}), is still unsolved even  for  quadratic
differential systems. The interest is in the global behavior of all solutions
in  the  whole plane and even at infinity (cf. \cite{Gonz}) and this for  a
whole family of systems, which is why this problem is so hard. The set \QS of
quadratic differential systems  depends on 12 parameters, the coefficients of
the two polynomials $p$ and $q$. On \QS acts the group
 of affine transformations  and of changes of scale on the time axis.
  So the  space  actually depends on five parameters. But even five is a
large number considering that we expect this class to yield
thousands of distinct phase portraits.  For  this reason people
have attempted to study particular classes  of  quadratic systems
and for some classes  a complete classification  of  phase
portraits with respect to topological equivalence was obtained
(quadratic systems  with a center \cite{Vlp}, \cite{DAS},
\cite{Pal_DS1}, quadratic Hamiltonian systems  \cite{Art_Lb1},
\cite{Cal_Vlp}, quadratic chordal systems \cite{Gsl}, quadratic
systems with a weak focus of third order \cite{Art_Lb2} and
\cite{Lib_DS}, etc.).

The  goal  in  most  of  these  articles  was  to obtain all topologically
distinct phase portraits for the particular class considered. Two systems ($S$)
and ($S'$) are topologically equivalent if there exists a homeomorphism $f:\
\R^2 \longrightarrow\ \R^2$ such that $f$ carries orbits to orbits preserving
(or reversing) their orientation. In most articles,  the classifications  were
done by using specific charts and normal forms for the systems in these charts
with respect  to parameters satisfying certain inequalities or equations.  The
results are not readily applicable  for systems  given in normal forms with
respect to other charts. This dependence on a specific normal form yields
results which are not geometric. Indeed, ever since Klein gave his famous
Erlangen   program,  we  are  used  to  calling  a  property geometric, if  it
is  invariant  under  the action of some group.  In  this sense, most of the
results obtained are not geometrical  since  they  are not independent of the
charts considered.

The orbit space of \QS under the group action is five-dimensional.
The global study of the class \QS  makes it necessary to use
normal forms in several charts. To obtain the global picture, we
need to be able to glue these charts and for this gluing process
invariants are helpful. Furthermore in work in progress on this
space, invariants are very helpful even in choosing the charts.

 Chart   independent  classification  results were obtained by K.S.
Sibirsky and his school (cf. \cite{Sib1}, \cite{Pop_Sib}, \cite{Cal_Vlp}) using
the algebraic invariant theory of differential equations developed by Sibirsky
and his disciples (cf. \cite{Sib2}, \cite{Sib3}, \cite{Vlp2}, \cite{Popa4}).
However most of the articles of the school of Sibirsky were published in
Russian, only some appeared in translations which partly explains why this
theory is rather unknown in the west.
 In these articles, invariants and comitants are
introduced  in  their  multi-index  tensorial
form, certain rather artificial
polynomial combinations of these
are  chosen  and classifications are given in terms of these
combinations.   In  the  end  these  classifications  remain
insufficiently  related  to  the  geometry  of  the systems.

We  need    much  simpler invariants,  simpler  than the configuration space of
Markus (cf.\cite{Mrk}),  possibly  even  integer  valued invariants which could
convey to us in simple terms properties of the global geometry  of  the
systems.  We  would  also  need  a way of computing  at  least some of these
simple invariants for any system in whatever chart it may be given to us.

In \cite{Pal_DS1}, \cite{Lib_DS} the authors gave topological
 classifications in terms of
the   global   geometry  of  the classes of systems considered.
These classifications  are affinely  invariant and they are
expressed in \cite{Pal_DS1} in  terms  of the geometry of
algebraic invariant curves of the systems considered  or in terms
of very simple  integer-valued invariants in \cite{Lib_DS} and
\cite{Pal_DS2} reflecting the geometry of the systems. There is
however a need to have an efficient way of effectively computing,
independent of charts these simpler integer-valued , geometric
invariants.

 In spite of their
awesome character, polynomial invariants  and  comitants are a very powerful
computational tool applicable to any canonical form and  they can be programmed
on a computer. There is thus a need to merge the geometric methods above
mentioned with the algebraic  invariant  approach. In this work we propose to
do just   that   for  the  specific  problem  of classifying topologically
quadratic systems  in the neighbourhood of  infinity.

     In \cite{Koj_Rn}  Kooij and Reyn  obtained all possible
local phase portraits  around  a single singular  point  at
infinity of an arbitrary quadratic  vector  field. In
\cite{Koj_Rn} they did not consider  the possible  ways  of
combining such singularities so as to obtain a topological
classification of quadratic systems in a neighborhood of the line
at infinity. In \cite{Nik_Vlp} I.~Nicolaev and N.~Vulpe obtained
such  a classification  in terms  of  algebraic  invariants  and
comitants and in \cite{Blt_Vlp1} the affine invariant
classification of quadratic system with respect to  the possible
distributions of the multiplicities of singularities at  infinity
was obtained by V.Baltag and N.Vulpe \cite{Blt_Vlp1}. These
classifications use the  technical language of algebraic invariant
theory developed by the school of Sibirsky
(\cite{Sib2},\cite{Vlp2},\cite{Bul_Tim}, etc).

 The goal of this work
is to combine the geometric approach in \cite{Pal_DS1},
\cite{Lib_DS} and \cite{Pal_DS2}  with the algebraic  invariant
approach  in \cite{Nik_Vlp} and \cite{Blt_Vlp1} for the
topological classification of quadratic systems in the
neighborhood of infinity.  We need  simplicity and clarity in the
geometric classification as well as applicability to any
particularly chosen chart.  In this article, which is based on
\cite{Dana_Vlp} we introduce the notions and prove the necessary
results which permit this  in as self-contained a way as possible.
We also point out that in the attempt to merge the geometric
invariants with the algebraic ones, the geometry led us to simpler
algebraic invariants than those in \cite{Nik_Vlp} yielding simpler
conditions in the classification Theorem \ref{th_2}.

This work could be applied along with an analogous one for finite
singularities, as an initial step, to the problem of classifying all quadratic
differential systems.

 The  article  is organized as follows:
In \S 2 we  consider  the two compactifications of real planar
polynomial systems and the foliations with singularities, real and
complex, on the real and complex projective planes, associated to these systems.

 In \S  3 we describe the purely geometric objects,
 i.e. the divisors
attached to the line at infinity, introduced in \cite{Pal_DS2}, which encode
  the multiplicities at infinity of
 the systems, and attach to these some integer-valued global
 affine invariants.

 In \S 4 we consider group actions on quadratic differential systems and
 define algebraic invariants and comitants with respect to these group
actions. We also give, using  a  comitant, canonical forms for these
differential systems according to their
 behavior at infinity.

In \S 5 we state and prove the classification theorem
(Theorem~\ref{th_1}) of the quadratic differential systems
according to their multiplicity divisors at infinity and for each
class we give the necessary and sufficient conditions in terms of
algebraic invariants and comitants with respect to the group
action. These conditions allow us to compute for any system and in
any chart the types of the multiplicity divisors associated to the
 system.

In \S 6 we introduce new classifying tools, among them the index divisor
encoding globally the topological indices of the singularities at
infinity of any polynomial differential system without a line of
singularities at infinity. We also introduce
a divisor encoding globally the number of local separatrices bounding
a hyperbolic sector of a singular point at infinity.

In \S 7 we  state and prove the topological classification theorem
(Theorem \ref{th_2}). This classification is expressed on one side
in terms of geometrical, affine integer-valued invariants, which
convey in simple terms the geometric and dynamic properties of the
systems according to their behavior in the vicinity of infinity;
on the other hand in terms of algebraic invariants and comitants.
In the end we are able to read for any system and in any chart,
its geometric and dynamical properties at
 infinity once these algebraic invariants and comitants are calculated.
 These calculations could  be done on a computer. A complete dictionary of
 integer-valued geometric invariants and polynomial invariants is given.

In  the  Appendix we list the invariants and comitants  used in \cite{Nik_Vlp}
and  which are needed for the proofs of the main results. These are also listed
for the purpose of comparison with the simpler algebraic  invariants and
comitants used in this article.  Highlighting the geometry of the systems via
the integer-valued invariants introduced, helped us to choose better algebraic
invariants and comitants than those in \cite{Nik_Vlp},
 closer to the geometry of the systems.

\vspace{-4mm}
\section{ The two  compactifications of real planar\\ polynomial
vector fields}
\vspace{-3mm}
A real planar polynomial system $(S)$ can be compactified on the sphere as
follows: Consider the $x,y$ plane as being the plane $Z = 1$ in the space
$\mathbb R^3$ with coordinates $X,Y,Z$. The central projection of the vector
field $p\partial/\partial x + q\partial /\partial y$ on the sphere of radius
one yields a diffeomorphic vector field on the upper hemisphere and also
another vector field on the lower
 hemisphere.
Poincar\'e indicated briefly in \cite{Po1} that one can construct  an analytic
vector field $\cal V$ on the whole sphere such that its restriction on the
upper hemisphere has the same phase curves as the one
 induced by
the phase curves of $(S)$ via the central projection. A complete
proof was given much later in \cite{Gonz}. The analytic vector
field $\cal V$ on the whole sphere obtained in this way  is called
the Poincar\'e field associated to the system $(S)$. The phase
curves of $\cal V$ coincide in each chart with phase curves
 induced by  planar
polynomial vector fields, in particular in the chart corresponding to $Z = 1$,
denoting the two coordinate axes $x, y$ corresponding to the $OX$ and $OY$
directions, they coincide with the phase curves induced by  $(S)$.  The two
planar polynomial vector fields $U, V$ associated to the
 charts for $X = 1$ (with local coordinates $(u,z)$) and for $Y = 1$
 (with local
coordinates $(v,w)$) and changes of coordinates $u = y/x$, $z =
1/x$,
 or $v = x/y,$ $w = 1/y$ are as
follows:
$$
 U\  \left\{\ba{ll} \displaystyle{\frac{du}{dt}} =
 &\displaystyle{C(1,u,z),}\\[2mm]
            \displaystyle{\frac{dz}{dt}}=&\displaystyle{ zP(1,u,z),}
      \ea \right.
\qquad
\mbox{and}\qquad
 V\ \left\{\ba{ll} \displaystyle{\frac{dv}{dt}} =
 &\displaystyle{C(v,1,w),}\\[2mm]
            \displaystyle{\frac{dw}{dt}}=&\displaystyle{-wQ(v,1,w),}
      \ea \right.
$$
where $P, Q$ and $C$ are defined further  below.

By the compactification of the planar polynomial vector field associated to
$(S)$ we understand the restriction ${\cal V}|_{\cal H'}$ (where by ${\cal H'}$
we understand the upper hemisphere ${\cal H}$ completed with the equator) of
the analytic vector field $\cal V$ on the sphere.  In this work we are
interested in the topological classification of $(S)$ on $\mathbb R^2$ (or
${\cal V}|_{\cal H}$) completed with its points "at infinity", i.e. on the
equator $S^1$ of $S^2$.  Since the vertical projection is a diffeomorphism of
${\cal H}'$ on the disk $\{ (x,y) | x^2 + y^2 \le1\}$  we can view the phase
portraits of our systems $(S)$ on this disk, called the Poincar\'e disk.

     We  shall  also  use  the compactifications (real or complex)
associated  to the foliations with singularities (real or complex) attached  to
a real polynomial system $(S)$ (cf. \cite{Cmch} or \cite{Dana_Sab}). These
foliations can  be  described  as  follows: For a real polynomial system  $(S)$
with  $n = \max(\deg(p), \deg(q))$  we associate to the two polynomials $p,
q\in \mathbb R[x,y]$ defining $(S)$, the homogeneous polynomials $P, Q$ in $X,
Y, Z$, of degree $n$ with real coefficients, defined as follows:
$$
  \ba{c} P(X,Y,Z) = Z^np(X/Z,Y/Z),\quad   Q(X,Y,Z) = Z^nq(X/Z,Y/Z).\\
   \ea
$$
     The real (respectively complex) foliations with singularities
associated  to  $(S)$  on  the real  projective plane
$\PP^2(\mathbb R)$  (respectively complex, $\PP^2(\mathbb C))$
are then described in homogeneous coordinates by the equation \bq
   && A(X,Y,Z)dX + B(X,Y,Z)dY + C(X,Y,Z)dZ = 0, \label{eq1}
\eq
where  $A = ZQ$, $B = - ZP$,\ $ C(X,Y,Z) = YP(X,Y,Z) - XQ(X,Y,Z) $ verify the
following equality\\
\vspace{-8mm}
\bq
   &&  A(X,Y,Z)X + B(X,Y,Z)Y + C(X,Y,Z)Z = 0  \label{eq2}
\eq in $\R[X,Y,Z]$.  (For more details see \cite{Cmch} or
\cite{Dana_Sab}).

 Our goal in this work is to give a topological classification,
in terms of  both geometric and algebraic invariants,
of  the  quadratic systems $(S)$  and their compactification on
$H'$ in the neighbourhood of the equator
in  the  closed  upper  hemisphere  $H'$  of  the  Poincar\'e  sphere.
Correspondingly  this  yields  a  topological classification of the
real  foliations,  in  the  neighbourhood  of the line at infinity
associated to the imbedding of the affine plane:
$$
                j:\ \AB^2(\mathbb R) = \mathbb R^2 \longrightarrow \PP^2(\mathbb R)
$$
where  $j(x,y) = [x:y:1]$. The line at infinity in this case is
therefore  $Z = 0$.
\vspace{-4mm}
\section{Divisors on the line at infinity encoding globally
the multiplicities of singularities}
\vspace{-3mm}

In this section we consider  real polynomial systems $(S)$ with  $ n  =
\max(\deg(p), \deg(q))$ and their associated foliations with singularities,
real or complex, defined in the previous section by the equation (\ref{eq1})
with coefficients $A,B,C$ verifying  (\ref{eq2}).
\bdf\label{df3_1}
For a system $(S)$  we call  divisor on the line at infinity,  a formal
expression of the form $D = \sum n(w)w$\  where $w$ is a point of the complex
line $Z=0$ of the complex projective plane, $n(w)$ is
 an integer and only a
finite number of the numbers $n(w)$ are not zero.  We call  degree
of the divisor $D$ the integer $\deg(D) = \sum n(w)$.  We call
support of the divisor $D$ the set $\Supp(D)$ of points $w$ such
that $n(w)\ne0$.
\edf
For systems $(S)$  two divisors on the line
at infinity were introduced in \cite{Pal_DS2}. These were applied
in \cite{Lib_DS} for classifying topologically the quadratic
systems with a weak focus of third order.
\bdf\label{df3_2}
 Assume that a system $(S)$ is such that $p(x,y)$ and $q(x,y)$ are
relatively prime over $\mathbb C$ and  that  $yp_n- xq_n$ is not identically
zero (i.e. $Z\nmid C$).  Here $p_n$  (respectively $q_n$) is the sum of terms
of degree $n$ of $p$ (respectively of $q$) in case at least one of them has a
non-zero coefficient and zero otherwise.

  The following divisor on the line at infinity is then well defined:
$$
       D_S(P,Q;Z) = \sum I_w(P,Q)w
$$
where the sum is taken for all points $w = [X:Y:0]$ on the line $Z = 0$ and
$I_w(P,Q)$ is the intersection number (or multiplicity of intersection) at $w$
(cf. \cite{Fult}) of the complex projective curves \quad $
                P(X,Y,Z) = 0\quad \hbox{and}\quad Q(X,Y,Z) = 0.
$
\edf
We thus have\qquad $
        \hbox{\rm Supp}(D_S(P,Q;Z)) = \{ w\in \{ Z = 0\} |
                  P(w) = 0 = Q(w)\}.
$

The above divisor is a purely geometric object which encodes the contribution
to the multiplicities of the singularities at infinity of the system $(S)$,
arising from singularities in the finite plane, i.e. how many singular points
in the finite plane could appear from those singularities at infinity in
(polynomial) perturbations  of $(S)$.

Let us list a number of integer-valued invariants which are attached to
this divisor.
\bnot
$$
\ba{c}
N_{\infty, f}(S) = \#\ \hbox{\rm Supp}(D_S(P,Q;Z));
\\
 \nu(S) = \max\{ I_w(P,Q)\, |\, w \in \hbox{\rm Supp}(D_S(P,Q;Z))\}; \\
  \hbox{for every}\ m \le \nu(S),\ s(m) = \#\ \{ w \in \{ Z = 0\}\, |\, I_w(P,Q) = m\}\ .
\ea
$$
\enot

Note that $N_{\infty,f}$ is the number of distinct infinite
singularities of $(S)$ which  could produce finite singular points
in a (polynomial) perturbation of $(S)$.

We also need another divisor on the line at infinity which was used in
\cite{Pal_DS2} and \cite{Lib_DS} and which is defined as follows:
\bdf Suppose $Z\nmid C$ and consider
$$
        D_S(C,Z) = \sum I_w(C,Z)w
$$
where the sum is taken for all points $w = [X:Y:0]$ on the line
$Z=0$ of the complex projective plane.
\edf

Clearly for quadratic differential systems $\deg(D_S(C,Z)) = 3$.
\bdf
 A point $w$ of the projective plane $\PP^2(\mathbb C)$ is said to be of multiplicity
 $(r,s)$  for a system $(S)$ if
 $$
                (r,s) = (I_w(P,Q),I_w(C,Z)).
$$
\edf

Following \cite{Pal_DS2} we  fuse the above two divisors on the line at infinity into just
one but  with values in the ring $\Z^2$:
\bdf
$$
        D_S = \sum\left(\ba{c}I_w(P,Q)\\ I_w(C,Z)\ea\right)w
$$
where $w$ belongs to the line $Z=0$ of the complex projective
plane.
\edf

The above defined divisor  describes the number of singularities which could
arise in a perturbation of $(S)$ from singularities at infinity of $(S)$ in
both the finite plane and  at infinity.\bdf
  We call type of the divisor $D_S(P,Q;Z)$  the set
$$
                \{ (s(m),\, m)\, |\, m \le \nu(S) \}.
$$
\edf

\brm
We observe that the types of $D_S(P,Q;Z)$ and of $D_S(C,Z)$ are affine invariants
since both $I_w(P,Q)$ and $I_w(C,Z)$ remain invariant under the action
of the affine group on systems $(S)$ (\cite{Pal_DS3}, \cite{DS}).
\erm
\bnot Let us introduce for planar systems $(S)$ the following notations:
$$
  \Delta_S = \deg D_S(P,Q;Z), \quad
  M_{C}=  \max \{I_w(C,Z)|\, w\in \hbox{\rm Supp}(D_S(C,Z))\}.
$$
\enot
\indent Consider a real quadratic differential system $(S)$:
\be
  \ba{ll}
\displaystyle  \frac {dx}{dt}&=p_0+ p_1(x,y)+\,p_2(x,y)\equiv p(x,y), \\[2mm]
\displaystyle  \frac {dy}{dt}&=q_0+ q_1(x,y)+\,q_2(x,y)\equiv q(x,y).
    \ea
    \label{2s1}
\ee
Suppose $\gcd(p,q)=\const$, where $p_i$ (respectively $q_i$) is the sum of
terms in $x$ and $y$ of degree $i$ of $p$ (respectively of $q$) in case at
least one such term has non-zero coefficient and zero otherwise. Recall that
\QS denotes the class of all real quadratic systems.

We want to list all possible divisors $D_S$ for quadratic systems $(S)$ and
characterize in terms of invariants and comitants the types of these divisors.
 This would make possible for any given system and in any chart
the computation of the type of its divisor $D_S$.
 To do this we need to construct invariants and
comitants with respect to group actions, which we do in the next section.
\section{Group actions on quadratic systems (\ref{2s1}) and invariants
and comitants with respect to these actions}\label{Sect_group}
\vspace{-4mm}
\subsection{Group actions on quadratic systems (\ref{2s1})}
More explicitly the systems (\ref{2s1}) can be written  in the form:
\beq
\frac{dx}{dt}&=& a_{00}+a_{10}x+a_{01}y+a_{20}x^2+2a_{11}xy+a_{02}y^2,\\
\frac{dy}{dt}&=& b_{00}+b_{10}x+b_{01}y+b_{20}x^2+2b_{11}xy+b_{02}y^2,
\eeq
and let   $ a =(a_{00},\ldots,b_{02}).$ Consider the ring
$\R[a_{00},a_{10},\ldots,a_{02},b_{00},b_{00},b_{10},\ldots,b_{02},\,x,y]$
which we shall denote  $\R[a,x,y]$.

On the set \QS of all quadratic differential systems (\ref{2s1}) acts the group
$\Aff(2,\mathbb R)$ of affine transformations on the plane. Indeed for every $
g\in \Aff(2,\mathbb R)$, $g:\ \mathbb R^{2}\longrightarrow \mathbb R^{2}$ we
have:
$$
  g:\ {\tilde x\choose \tilde y} =
      M\,{x\choose y} +B;\qquad
  g^{-1}:\ {x\choose y} =
      M^{-1}\, {\tilde x\choose \tilde y} -M^{-1}B.
$$
where $M=|| M_{ij} || $ is a $2\times 2$ nonsingular matrix and $B$ is a
$2\times 1$ matrix over $\mathbb R$. For every $ S\in \QS$ we can form its
transformed system $\tilde S=g S$:
$$
 \qquad\qquad \frac{d\tilde x}{dt} =\tilde p(\tilde x,\tilde y),\qquad\quad
  \frac{d\tilde y}{dt} =\tilde q(\tilde x,\tilde y),\eqno(\tilde S)
$$
where
$$
  \left(\ba{c} \tilde p(\tilde x,\tilde y)\\ \tilde q(\tilde x,\tilde y)\ea\right) =
      M\left(\ba{c} (p\,{\mbox{\footnotesize$\circ$}}\, {g^{-1}})(\tilde x,\tilde y)\\
                   (q\,{\mbox{\footnotesize$\circ$}}\, {g^{-1}})(\tilde x,\tilde y) \ea\right).
$$
The map
\beq
   \Aff(2,\mathbb R)\times \QS\ &\longrightarrow &\
   \QS \\
   (g,\quad S)\ &\longrightarrow &\ \tilde S=g S
\eeq
verifies the axioms for a left group action. For every  subgroup $G\subseteq
\Aff(2,\mathbb R)$ we have an  induced action of $G$ on \QS. We can identify
the set \QS of system (\ref{2s1}) with a subset of $\mathbb R^{12}$
 via the embedding  $\QS \hookrightarrow\R^{12}$ which associates to
 each system (\ref{2s1}) the 12-tuple $(a_{00},\ldots,b_{02})$ of its coefficients.

 On systems $(S)$ such that
$\max(\deg(p),\deg(q))\le2$ we consider the  action of the group
$\Aff(2,\R)$ which   yields an action of this group on $\mathbb
R^{12}$. For every $ g\in \Aff(2,\R)$ let $r_g:\
\R^{12}\longrightarrow \R^{12}$ be the map which corresponds to
$g$ via this action. We know \mbox{(cf. \cite {Sib1})} that $r_g$
is linear and that the map $r:\ \Aff(2,\R)\longrightarrow
GL(12,\R)$ thus obtained is a group homomorphism. For every
subgroup $G$ of $\Aff(2,\R)$, $r$ induces a representation of $G$
onto a subgroup  $\cal G$ of $GL(12,\R)$.
\subsection{Invariants and comitants associated to the group actions}
\vspace{-3mm} \bdf \label{def:comit} A polynomial $U( a\,,x,y)\in
\mathbb R[a,x,y]$ is called a comitant  of systems (\ref{2s1})
with respect to a subgroup  $G$ of $\Aff(2,\mathbb R)$, if there
exists $\chi\in \Z$ such that for every \mbox{$(g,\,\ab )\in
G\times\mathbb R^{12} $} and for every  $(x,y)\in \mathbb R^2$ the
following relation holds:
$$
U(r_g(\ab) ,\ g(x,y)\,)\equiv\ (\det\,g)^{-\chi}\, U( \ab ,x,y),
$$
where $\det g=\det M$.
If the polynomial $U$ does not explicitly depend  on
$x$ and $y$ then it is called invariant.
The number $ \chi\in \Z $  is called
the   {\sl weight\,} of the comitant $ U( a ,x,y)$.
If $G=GL(2,\mathbb R)$ (or $G=Aff(2,\mathbb R)$\,) then the comitant
$U( a ,x,y)$ of systems (\ref{2s1}) is called $GL$-comitant (respectively,
 affine comitant).
\edf
\bdf\label{def_G}
A subset $X\subset \mathbb R^{12}$ will be called $G$-invariant,\ if \ for
every $ g\in G$ we have $r_g(X)\subseteq X$.
\edf
As it can easily be verified,
the following  polynomials are $GL$-comitants
 of system (\ref{2s1}):
\be\label{com1}
\bal
   {C}_i(a,x,y)\ =\ & y p_i(x,y)-x q_i(x,y),\  i=0,1,2;\\
   M(a,x,y)\ =\ &  2\,\Hess\big(C_2(a,x,y)\big);\\
    \eta(a)\ =\ & \Discrim\big(C_2(a,x,y)\big);\\
   K(a,x,y)\ =\ &  \Jacob\big(p_2(x,y),q_2(x,y)\big);\\
   \mu_0(a)\ =\ & \Res_x(p_2,q_2)/y^4 \  = \ \Discrim\big(K(a,x,y)\big)/16;\\
   H(a,x,y)\  =\ &  -\left.\Discrim(\alpha p_2(x,y)+\beta q_2(x,y))
   \right|_{\{\alpha=y,\, \beta=-x\}};\\
   L(a,x,y)\ =\ & 2K-4H-M;\\
   K_1(a,x,y)\ =\ & p_1(x,y)q_2(x,y)-p_2(x,y)q_1(x,y).\\
\eal
\ee

Let $T(2,\R)$ be the subgroup of $\Aff(2,\R)$ formed by translations. Consider
the linear representation of $T(2,\R)$ into its corresponding subgroup ${\cal
T}\subset GL(12,\mathbb R)$, i.e. for every $ \tau\in T(2,\R)$,\ $\tau:\
x=\tilde x+\alpha, y=\tilde y+\beta$ we consider as above $r_\tau:\
\R^{12}\longrightarrow \R^{12}$.
\bdf
A $GL$-comitant $U( a ,x,y)$
of systems (\ref{2s1}) is called a $T$-comitant if for every
$(\tau,\, \mbox{\boldmath $a$} )\in T(2,\mathbb R)\times \mathbb R^{12}$ and
 for every $(\tilde x,\tilde y)\in \mathbb R^2$
the  relation \ $ U(r_\tau\cdot  \mbox{\boldmath $a$},\ \tilde x,\,\tilde y)\,
=\, U( \mbox{\boldmath $a$} ,\ \tilde x,\,\tilde y) $ holds.
\edf
Let
$$
U_i( a ,x,y)=\sum_{j=0}^{d_i} U_{ij}( a )x^{d_i-j}y^j,\quad
i=1,\ldots,s
$$
be a set of $GL$-comitants of systems (\ref{2s1}) where $d_i$ denotes the degree
of the binary form $U_i( a ,x,y)$ in $x$ and $y$ with coefficients in
$\mathbb R[a]$ where $\mathbb R[a]=\mathbb R[a_{00},\ldots,b_{02}]$.
We denote  by
$$
{\cal U}=\left\{\,U_{ij}( a )\in \mathbb R[a]\ |\ i=1,\ldots,s,\
j=0,1,\ldots,d_i \,\right\},
$$
 the set of the coefficients in $\R[a]$ of the $GL$-comitants $U_i( a ,x,y)$,
$ i=1,\ldots,s$, and by $V(\cal U)$ its associated algebraic set:
$$
  V({\cal U})=\left\{\,\mbox{\boldmath $a$}\in \mathbb R^{12}\ |\
  U_{ij}(\mbox{\boldmath $a$})=0 \ \forall\ U_{ij}( a )\in \cal U\,\right\}.
$$

\bdf\label{def:TC-comit}
 A $GL$-comitant $U( a ,x,y)$
of systems (\ref{2s1}) is called a conditional $\ T$-comitant (or $CT$-comitant)
 modulo $\left<U_1,U_2,...,U_s\right>$ if the following two conditions are
 satisfied:

 (i) the algebraic subset $V({\cal U})\subset \mathbb R^{12}$ is affinely  invariant
(see Definition \ref{def_G});

 (ii) for every ($\tau,\  \mbox{\boldmath $a$} )\in T(2,\mathbb R)\times V(\cal U)$\
 we have
$
U(r_\tau\cdot  \mbox{\boldmath $a$} ,\ \tilde x,\,\tilde y)=
 U( \mbox{\boldmath $a$} ,\ \tilde x,\,\tilde y)\
\mbox{in}\ \mathbb R[\tilde x,\,\tilde y]. $ \edf In other words,
a $CT$-comitant $U( a ,x,y)$ modulo $\left<U_1,U_2,...,U_s\right>$
is a $T$-comitant on the algebraic subset \mbox{$V({\cal
U})\subset \mathbb R^{12}$}.

\bpr\label{prop:p2,q2} Let $S\in \QS$ and let $\ab\in \R^{12}$ be its 12-tuple of
coefficients.  The common points of $P=0$ and $Q=0$ on the line $Z=0$ are given
by the common linear factors over $\C$ of $p_2$ and $q_2$. This yields the
geometrical meaning of the comitants  $\mu_0$, $K$ and $H$:
\beq
\gcd (p_2(x,y),q_2(x,y)) =\left\{\begin{array}{lcl}
           \const & iff & \ \mu_0(\ab)\ne0;\\
           bx+cy& iff & \ \mu_0=0,\ K( \ab ,x,y)\ne0;\\[1mm]
          (bx+cy)(dx+ey) & iff &\rule{0pt}{3.5ex}  \left\{\ba{l}\mu_0(\ab)=0, K( \ab ,x,y)=0\\
                                                \hbox{and}\ \ H( \ab ,x,y)\ne0;\\
                                        \ea\right. \\[4mm]
          (bx+cy)^2 & iff &  \rule{0pt}{4ex} \left\{\ba{l}\mu_0=0, K( \ab ,x,y)=0,\\
                                               \hbox{and}\ \  H( \ab ,x,y)=0;\\
                                        \ea\right. \\
          \end{array}\right.
\eeq
where $bx+cy,dx+ey\in \mathbb C[x,y]$ are some linear forms and $be-cd\ne0$.
\epr
\bdf
 The polynomial $U( a ,x,y)\in \mathbb R[ a ,x,y]$
 has  well determined sign on \mbox{$V\subset\mathbb R^{12}$} with respect to $x,\,y$
 if for every fixed $\ab\in V$, the sign of the polynomial function  $U(\ab,x,y)$
   on $\R^2$ is constant where this function is not zero.
\edf
\bobs We draw the attention to the fact, that if a $CT$-comitant
$U( a ,x,y)$ of even weight is a binary form in $x$, $y$, of even
degree in the coefficients  of (\ref{2s1}) and has well determined
sign on some affine invariant algebraic subset $V(\cal U)$ then
this property is conserved by any affine transformation and the
sign is conserved. \eobs
\vspace{-4mm}
\subsection{Canonical forms of planar quadratic systems
in the\\ neighbourhood of infinity}
\vspace{-3mm}
\blm\label{lm2}
For a system (\ref{2s1}) with $C_2( \ab ,x,y)\not\equiv 0$  the
divisor $D_S(C,Z)$ is well defined and its type is determined by
the corresponding conditions indicated in Table 1, where we write
$q_1^c+q_2^c+q_3$ if two of the points, i.e. $q_1^c, q_2^c$, are
complex but not real.
        Moreover, for each type of the divisor $D_S(C,Z)$ given
by Table 1 the quadratic systems (\ref{2s1}) can be brought via a
linear transformation to one of the following  canonical systems
$(\SSS_{I})-(\SSS_{IV})$ corresponding to their behavior at
infi\-ni\-ty.
\elm
\begin{table}[!htb]
\begin{center}
\begin{tabular}{|c|c|c|c|}
\multicolumn{4}{r}{\bf Table 1}\\[2mm]
\hline
  \raisebox{-0.7em}[0pt][0pt]{$M_{C}$ }  & \raisebox{-0.7em}[0pt][0pt]{Type of $D_S(C,Z)$}
      & Necessary and sufficient &  Notation for  \\[-1mm]
         & & conditions on the comitants & the conditions\\
 \hline\hline
 \rule{0pt}{1.2em}
 \raisebox{-1.0em}[0pt][0pt]{$1$} & $q_1+q_2+q_3 $ &  $\eta>0 $ & $({\cal I}_1)$ \\[1mm]
\cline{2-4}
 \rule{0pt}{1.2em}   & $q_1^c+q_2^c+q_3 $ &  $\eta<0$ & $({\cal I}_2)$\\[1mm]
\hline
 \rule{0pt}{1.2em}  $2$ & $2q_1+q_2 $ &  $\eta=0,\quad M\ne0$ & $({\cal I}_3)$\\[1mm]
\hline
 \rule{0pt}{1.2em} $3$ & $3q $ &  $ M=0$ &$({\cal I}_4)$ \\[1mm]
\hline
\end{tabular}
\end{center}
$$
\bal
&\qquad\left\{\ba{rcl}
 \displaystyle \frac{dx}{dt}&=&k+cx+dy+gx^2+(h-1)xy,\\[2mm]
 \displaystyle \frac{dy}{dt}&=& l+ex+fy+(g-1)xy+hy^2;
\ea\right. \qquad &\qquad (\SSS_I)\\[3mm]
&\qquad\left\{\ba{rcl}
 \displaystyle \frac{dx}{dt}&=&k+cx+dy+gx^2+(h+1)xy,\\[2mm]
 \displaystyle \frac{dy}{dt}&=& l+ex+fy-x^2+gxy+hy^2;
\ea\right.& \hspace{2cm}(\SSS_{I\!I})\\[3mm]
&\qquad\left\{\ba{rcl}
 \displaystyle \frac{dx}{dt}&=&k+cx+dy+gx^2+hxy,\\[2mm]
 \displaystyle \frac{dy}{dt}&=& l+ex+fy+(g-1)xy+hy^2;
\ea\right.&\qquad (\SSS_{I\!I\!I})\\[3mm]
&\qquad\left\{\ba{rcl}
 \displaystyle \frac{dx}{dt}&=&k+cx+dy+gx^2+hxy,\\[2mm]
 \displaystyle \frac{dy}{dt}&=& l+ex+fy-x^2+gxy+hy^2,
\ea\right.&\qquad (\SSS_{I\!V})\\[3mm]
\eal
\vspace{-5mm}
$$
\end{table}
\BProof The Table 1 follows easily from the definitions of  $\eta(a)$ and  $
M(a,x,y)$ in \eqref{com1}. Let us consider the $GL$-comitant
${C_2}( a ,x,y)\not\equiv0$ simply as a cubic binary form in $x$
and $y$. For every $ \ab\in\R^{12}$ the binary form ${C_2}( \ab
,x,y)$ can be reduced to one of the  canonical forms given below,
by a linear transformation, i.e. there exist $g\in GL(2,\R)$:
$g(x,y) = (u,v)$ such that the transformed binary form $g{C_2}(
\ab,x,y)=C_2(\ab,g^{-1}(u,v))$ is one of the following
\be\label{can_forms}
I.\ xy(x-y); \qquad II.\ x(x^2+y^2); \qquad III.\ x^2y; \qquad IV.\ x^3,
\ee
which correspond to the types of the divisor $D_S(C,Z)$ indicated in Table
 1. On the other hand, according to the Definition
\ref{def:comit} of the $GL$-comitant, for ${C_2}(a ,x,y)$ whose weight
$\chi=-1$, we have for $g\in GL(2,\R)$
$$
  C_2(r_g(\ab),\, g(x,y))= \det(g)\, C_2( \ab,\, x,y).
$$
Using $g(x,y) = (u,v)$ we obtain
$$
C_2(r_g(\ab),\, u,v)=\lambda C_2( \ab ,\, g^{-1}(u,v)), \quad \lambda \in \R,
$$
where we may consider $\lambda=1$ by rescaling :  $u=u_1/\lambda, \
v=v_1/\lambda$.

Thus, recalling that
$$
p_2(x,y)=a_{20}x^2+2a_{11}x,y+a_{02}y^2,\qquad
q_2(x,y)=b_{20}x^2+2b_{11}x,y+b_{02}y^2,
$$
  for the first canonical form in (\ref{can_forms}) we have
$$
C_2(\ab,x,y)=-b_{20}x^3+(a_{20}-2b_{11})x^2y+(2a_{11}-b_{02})xy^2+a_{02}y^3=
xy(x-y).
$$
Identifying the coefficients of the above identity we get the canonical form
$(\SSS_{I})$: Analogously for the cases $II,\ III$ and $IV$ we obtain the
canonical form $(\SSS_{I\!I})$, $(\SSS_{I\!I\!I})$ and $(\SSS_{I\!V})$
associated to the respective polynomials in (\ref{can_forms}). \EProof

\section{Classification of the quadratic systems according to the types of
the multiplicity divisor ${\cal D}_S$} \label{Sect_infty}
A specific type of a divisor  $ D_S$ yields a class of quadratic systems
(\ref{2s1}). We want to list all possible types of the divisors $ D_S$ and for
each specific type to determine the subset of \QS where $ D_S$ has this type.
We want to give this subset in terms of algebraic invariants and comitants so
as to be able to check these conditions for every
 system (\ref{2s1}) in any chart.

In order to construct other necessary invariant polynomials
let us consider the differential operator
${\cal L}= x\cdot \LL_2 -y\cdot\LL_1$
acting on $\mathbb R[a,x,y]$ constructed in \cite{Blt_Vlp2},
where
\beq
  && \LL_1= 2a_{00}\frac{\partial}{\partial a_{10}} +
            a_{10}\frac{\partial}{\partial a_{20}} +
    \frac{1}{2}a_{01}\frac{\partial}{\partial a_{11}} +
            2b_{00}\frac{\partial}{\partial b_{10}} +
            b_{10}\frac{\partial}{\partial b_{20}} +
     \frac{1}{2}b_{01}\frac{\partial}{\partial b_{11}};\\
  && \LL_2= 2a_{00}\frac{\partial}{\partial a_{01}} +
            a_{01}\frac{\partial}{\partial a_{02}} +
     \frac{1}{2}a_{10}\frac{\partial}{\partial a_{11}} +
            2b_{00}\frac{\partial}{\partial b_{01}} +
            b_{01}\frac{\partial}{\partial b_{02}} +
     \frac{1}{2}b_{10}\frac{\partial}{\partial b_{11}}
\eeq
as well as the classical  differential operator
$(f,\varphi)^{(2)}$ acting on $\mathbb R[a,x,y]$ which is called
{\sl transvectant} of the second index (see, for example,
\cite{Gr_Yng,Olver}):
\beq
(f,\varphi)^{(2)}= \frac{\partial^2 f}{\partial x^2}
             \frac{\partial^2 \varphi}{\partial y^2}
           -2\frac{\partial^2 f}{\partial x\partial y}
             \frac{\partial^2 \varphi}{\partial x\partial y}
            +\frac{\partial^2 f}{\partial y^2}
             \frac{\partial^2 \varphi}{\partial x^2}.
\eeq
Here $f(x,y)$ and $\varphi(x,y)$ are polynomials in $x$ and $y$.

In \cite{Blt_Vlp2} it is shown that if a polynomial $U\in \mathbb
R[a,x,y]$ is a comitant of system  (\ref{2s1}) with respect to the
group $GL(2,\mathbb R)$ then ${\cal L}(U)$ is also a
$GL$-comitant.  The same is true for the operator transvectant  of
two comitants, $f$ and $\varphi$.

So, by using these operators and the $GL$-comitants $\mu_0(a),
M(a,x,y)$ and $K(a,x,y)$ we shall construct the following
polynomials:
\be\label{com2}
    \mu_i(a,x,y) =\frac{1}{i!} {\cal L}^{(i)}(\mu_0), \ i=1,..,4,\ \quad
     \kappa(a)=(M,K)^{(2)}, \quad   \kappa_1(a)=(M,C_1)^{(2)},
\ee
where ${\cal L}^{(i)}(\mu_0)={\cal L}({\cal L}^{(i-1)}(\mu_0))$.

These polynomials are in fact comitants of system (\ref{2s1}) with respect to the
group $GL(2,\mathbb R)$.

To reveal the geometrical meaning of the comitants $\mu_i(a, x, y),\
i=0,1,\ldots,4 $ we use  the following resultants whose calculation yield:
\bq
Res_{{}_X}\left(P, Q\right)&=& \mu_0 Y^4+\mu_{10}Y^3Z+
        \mu_{20}Y^2Z^2+\mu_{30}YZ^3+\mu_{40} Z^4;  \label{2s2}\\
Res_{{}_Y}\left(P, Q\right)&=& \mu_0 X^4+\mu_{01}X^3Z+
      \mu_{02}X^2Z^2+ \mu_{03}XZ^3+\mu_{04} Z^4,\quad \label{2s3}
\eq
where $\mu_{ij}=\mu_{ij}(a)\in \mathbb R[a_{00},\ldots, b_{02}]$.

On the other hand for $\mu_i$,\ $i= 0,1,\ldots,4$ from (\ref{com2}) we have
\beq
\begin{array}{lcc}
\mu_0(a)&=&\mu_0;\\[1mm]
\mu_1(a,x,y)&=&\mu_{10}x+\mu_{01}y;\\[1mm]
\mu_2(a,x,y)&=&\mu_{20}x^2+\mu_{11}xy+\mu_{02}y^2;\\[1mm]
\mu_3(a,x,y)&=&\mu_{30}x^3+\mu_{21}x^2y+\mu_{12}xy^2+\mu_{03}y^3;\\[1mm]
\mu_4(a,x,y)&=&\mu_{40}x^4+\mu_{31}x^3y+\mu_{22}x^2y^2+\mu_{13}xy^3+\mu_{04}y^4.
\end{array}
\eeq
We observe that  the leading coefficients of the comitants $\mu_i$,\ $i=
0,1,\ldots,4$ with  respect to $x$ (respectively $y$) are the corresponding
coefficients in (\ref{2s2})  (respectively (\ref{2s3})).

We draw the attention to the fact, that if the comitant $\mu_i(a,x,y)$\
$(i=0,1,\ldots,4)$ is not equal to zero then we  may assume that its leading
coefficients are both non zero, as this can be obtained by applying a rotation
of the phase plane of the system (\ref{2s1}). From here and (\ref{2s2}),
(\ref{2s3}) and the above values of $\mu_i$,\ $i= 0,1,\ldots,4$ we have:
\blm\label{lm_1}
The system $P(X,Y,Z)=Q(X,Y,Z)=0$  possesses  $m\ (=\Delta_S)$ ($1\le m \le4$)
solutions $[X_i:Y_i:Z_i]$ with $Z_i=0$ $(i=1,\ldots,m$) (considered with
multiplicities) if and only if for every $ i\in\left\{0,1,\ldots,m-1\right\}$
we have $\mu_i(a,x,y)=0$ in $\mathbb R[a,x,y]$ and $\mu_m(a,x,y)\ne0$.
\elm
\brm\label{rm1} It can easily be checked that the following identity holds\\
\vspace{-2mm}
$$
\vspace{-2mm}
 \mu_4(a,X,Y)\ =\ Res_Z\left(P(X,Y,Z), Q(X,Y,Z)\right).
$$
Hence, clearly for any solution $[X_0:Y_0:Z_0]$ (including those with  $Z_0=0$)
of the system of equations $P(X,Y,Z)=Q(X,Y,Z)=0$, the following relation is
satisfied: $ \mu_4(a,X_0,Y_0)=0. $
\erm
\vspace{-3mm}

We give below our theorem of classification of the types of all divisors $D_S$
occurring in quadratic systems and associate to each type the necessary and
sufficient conditions
 in terms of algebraic
invariants  and  comitants. The computation of these invariants and comitants
 can be
programmed  using symbolic manipulations and implemented
on  computers.  Thus  for any specific system (\ref{2s1}) we can calculate
explicitly its divisor type in whatever chart (\ref{2s1}) is given.
\bth\label{th_1}
    We consider here the family {\bf QS}$_{\bf ess}$  of all systems $(S)$ in \QS which
are essentially quadratic,  i.e. $\gcd(P,Q)=1$ and $Z\nmid C$.
 All possible values which could be
taken by $\Delta_S$ for such systems (\ref{2s1}) are as listed in
the first column of \mbox{Table 2}. For each value of $\Delta_S$,
all possibilities we have for $M_{C}$, are  listed in the second
column. For each combination $(\Delta_S, M_{C})$ all the
possibilities we have for the form of $D_S$ are those indicated in
the third column. For a specified $(\Delta_S, M_{C})$, the
necessary and sufficient conditions to have the form of $D_S$ as
indicated in the third column are those indicated in the
corresponding fourth column. (We recall that ${\cal I}_j$ are the
conditions indicated in Table 1. In the last column of Table 2 we
denote by $\Sigma_i$ the class of all quadratic systems which
possess $(\Delta_S, M_{C}, D_S)$ as indicated in the first three
columns).
\eth
\begin{table}
\begin{tabular}{|c|c|c|l|c|}
\multicolumn{5}{r}{\bf Table 2 }\\
\multicolumn{5}{r}{}\\
\hline
 \raisebox{-0.8em}[0pt][0pt]{$\Delta_S $}  &  \raisebox{-0.8em}[0pt][0pt]{$\!M_{C }\!$ }&
 \raisebox{-0.8em}[0pt][0pt]{Value of ${ D}_S$ }&\hfil Necessary and sufficient\hfill &
 \raisebox{-0.8em}[0pt][0pt]{$\Sigma_i$}\\
 & &  & \hfil conditions on the comitants \hfill & \\
 \hline\hline
 \rule{0pt}{2.8ex}  &    \raisebox{-0.7em}[0pt][0pt]{1} &$ {0\choose 1}p+{0\choose 1}q+{0\choose 1}r$
                                           & $ \mu_0\ne0,\ ({\cal I}_1)$& $\Sigma_{1}$\\[0.7mm]
\cline{3-5} \rule{0pt}{2.8ex}   &      &$ {0\choose 1}p+{0\choose
1}q^c+{0\choose 1}r^c$
                                           & $ \mu_0\ne0,\ ({\cal I}_2)$&
                                           $\Sigma_{2}$\\[0.7mm]
\cline{2-5} \rule{0pt}{2.8ex}\raisebox{0.7em}[0pt][0pt]{0}    & 2
&$ {0\choose 1}p+{0\choose 2}q $
                                           & $ \mu_0\ne0,\ ({\cal I}_3)$& $\Sigma_{3}$\\[0.7mm]
\cline{2-5}
\rule{0pt}{2.8ex}   &    3 &$ {0\choose 3}p $   & $ \mu_0\ne0,\ ({\cal I}_4)$& $\Sigma_{4}$\\[0.7mm]
\hline \rule{0pt}{2.8ex}   &    \raisebox{-0.7em}[0pt][0pt]{1} &$
{1\choose 1}p+{0\choose 1}q+{0\choose 1}r$
                                  & $ \mu_0=0,\ \mu_1\ne0,\ ({\cal I}_1)$& $\Sigma_{5}$\\[0.7mm]
\cline{3-5}
 \rule{0pt}{2.8ex}  &     &$ {1\choose 1}p+{0\choose 1}q^c+{0\choose 1}r^c$
                                   & $ \mu_0=0,\ \mu_1\ne0,\ ({\cal I}_2)$& $\Sigma_{6}$\\[0.7mm]
\cline{2-5} \rule{0pt}{2.8ex} 1 &
\raisebox{-0.7em}[0pt][0pt]{2} &$ {1\choose 1}p+{0\choose 2}q $
                          & $ \mu_0=0,\ \mu_1\ne0,\ \kappa\ne0,\ ({\cal I}_3)$& $\Sigma_{7}$\\[0.7mm]
\cline{3-5} \rule{0pt}{2.8ex}   &      &$ {0\choose 1}p+{1\choose 2}q $
                            & $ \mu_0=0,\ \mu_1\ne0,\ \kappa=0,\ ({\cal I}_3)$& $\Sigma_{8}$\\[0.7mm]
\cline{2-5}
\rule{0pt}{2.8ex}   &    3 &$ {1\choose 3}p $   & $ \mu_0=0, \ \mu_1\ne0,\ ({\cal I}_4)$& $\Sigma_{9}$\\[0.7mm]
\hline \rule{0pt}{2.8ex}   &      &$ {2\choose 1}p+{0\choose 1}q+{0\choose 1}r$
                   & $ \mu_{0,1}=0, \ \mu_2\ne0,\ \kappa\ne0,\ ({\cal I}_1)$& $\Sigma_{10}$\\[0.7mm]
\cline{3-5} \rule{0pt}{2.8ex}   &   \raisebox{-0.7em}[0pt][0pt]{1}
&$ {1\choose 1}p+{1\choose 1}q+{0\choose 1}r$
                   & $ \mu_{0,1}=0,\ \mu_2\ne0,\ \kappa=0,\ ({\cal I}_1)$& $\Sigma_{11}$\\[0.7mm]
\cline{3-5} \rule{0pt}{2.8ex}   &     &$ {2\choose 1}p+{0\choose
1}q^c+{0\choose 1}r^c$
                   & $ \mu_{0,1}=0,\ \mu_2\ne0,\ \kappa\ne0,\ ({\cal I}_2)$& $\Sigma_{12}$\\[0.7mm]
\cline{3-5} \rule{0pt}{2.8ex} \raisebox{-0.7em}[0pt][0pt]{2}  & &$
{0\choose 1}p+{1\choose 1}q^c+{1\choose 1}r^c$
                   & $ \mu_{0,1}=0,\ \mu_2\ne0,\ \kappa=0, \ ({\cal I}_2)$& $\Sigma_{13}$\\[0.7mm]
\cline{2-5} \rule{0pt}{2.8ex}   &     &$ {2\choose 1}p+{0\choose 2}q $
                   & $ \mu_{0,1}=0,\ \mu_2\ne0,\ \kappa\ne0,\ ({\cal I}_3)$& $\Sigma_{14}$\\[0.7mm]
\cline{3-5} \rule{0pt}{2.8ex}    &   2   &$ {1\choose 1}p+{1\choose 2}q $
                   & $ \mu_{0,1}=0,\ \mu_2\ne0,\ \kappa=0,\  L=0,\ ({\cal I}_3)$& $\Sigma_{15}$\\[0.7mm]
\cline{3-5} \rule{0pt}{2.8ex}   &    &$ {0\choose 1}p+{2\choose 2}q $
                   & $ \mu_{0,1}=0,\ \mu_2\ne0,\ \kappa=0,\ L\ne0,\ ({\cal I}_3)$& $\Sigma_{16}$\\[0.7mm]
\cline{2-5}
\rule{0pt}{2.8ex}   &   3 &$ {2\choose 3}p $   & $ \mu_{0,1}=0,\ \mu_2\ne0,\ ({\cal I}_4)$& $\Sigma_{17}$\\[0.7mm]
\hline \rule{0pt}{2.8ex}    &    &$ {3\choose 1}p+{0\choose 1}q+{0\choose 1}r$
                  & $ \mu_{0,1,2}=0,\ \mu_3\ne0,\ \kappa\ne0,\ ({\cal I}_1)$& $\Sigma_{18}$\\[0.7mm]
\cline{3-5}
 \rule{0pt}{2.8ex}  &   1  &$ {2\choose 1}p+{1\choose 1}q+{0\choose 1}r$
                   & $ \mu_{0,1,2}=0,\ \mu_3\ne0,\ \kappa=0,\ ({\cal I}_1)$& $\Sigma_{19}$\\[0.7mm]
\cline{3-5} \rule{0pt}{2.8ex}   &     &$ {3\choose 1}p+{0\choose
1}q^c+{0\choose 1}r^c$
                   & $ \mu_{0,1,2}=0,\ \mu_3\ne0, \ ({\cal I}_2)$& $\Sigma_{20}$\\[0.7mm]
\cline{2-5} \rule{0pt}{2.8ex} \raisebox{-0.7em}[0pt][0pt]{3} & &$
{3\choose 1}p+{0\choose 2}q $
             & $ \mu_{0,1,2}=0,\ \mu_3\ne0,\ \kappa\ne0,\ ({\cal I}_3)$& $\Sigma_{21}$\\[0.7mm]
\cline{3-5} \rule{0pt}{2.8ex}   &   \raisebox{-0.7em}[0pt][0pt]{2}
&$ {2\choose 1}p+{1\choose 2}q$
          & $ \mu_{0,1,2}=0, \mu_3\ne0, \kappa= L=0, \kappa_1\ne0,  ({\cal I}_3)$& $\Sigma_{22}$\\[0.7mm]
\cline{3-5} \rule{0pt}{2.8ex}   &   &$ {1\choose 1}p+{2\choose 2}q$
                   & $ \mu_{0,1,2}=0, \mu_3\ne0, \kappa= L=0,  \kappa_1=0,  ({\cal I}_3)$& $\Sigma_{23}$\\[0.7mm]
\cline{3-5}
 \rule{0pt}{2.8ex}  &    &$ {0\choose 1}p+{3\choose 2}q $
             & $ \mu_{0,1,2}=0,\ \mu_3\ne0,\ \kappa=0,\ L\ne0,\ ({\cal I}_3)$& $\Sigma_{24}$\\[0.7mm]
\cline{2-5}
\rule{0pt}{2.8ex}  &   3 &$ {3\choose 3}p $   & $ \mu_{0,1,2}=0,\ \mu_3\ne0,\ ({\cal I}_4)$& $\Sigma_{25}$\\[0.7mm]
\hline
\end{tabular}
\end{table}

\begin{table}
\begin{tabular}{|c|c|c|l|c|}
\multicolumn{5}{r}{{\bf Table 2}\ (continued)}\\
\multicolumn{5}{r}{}\\
\hline
 \raisebox{-0.8em}[0pt][0pt]{$\Delta_S $}  &  \raisebox{-0.8em}[0pt][0pt]{$\!M_{C }\!$ }&
 \raisebox{-0.8em}[0pt][0pt]{Value of ${ D}_S$ }&\hfil Necessary and sufficient\hfill &
 \raisebox{-0.8em}[0pt][0pt]{$\Sigma_i$}\\
 & &  & \hfil conditions on the comitants \hfill & \\
 \hline\hline
\rule{0pt}{2.8ex}    &    &$ {4\choose 1}p+{0\choose 1}q+{0\choose 1}r$
                   & $ \mu_{0,1,2,3}\!=\!0,\ \mu_4\ne0,\ \kappa\ne0,\ ({\cal I}_1)$& $\Sigma_{26}$\\[0.7mm]
\cline{3-5} \rule{0pt}{2.8ex}   &    &$ {3\choose 1}p+{1\choose 1}q+{0\choose
1}r$
                   & $ \mu_{0,1,2,3}\!=\!0,\ \mu_4\ne0,\ \kappa=0,\ K_1\ne0,\ ({\cal I}_1)$& $\Sigma_{27}$\\[0.7mm]
\cline{3-5} \rule{0pt}{2.8ex}   & 1    &$ {2\choose 1}p+{2\choose 1}q+{0\choose
1}r$
                   & $ \mu_{0,1,2,3}\!=\!0,\ \mu_4\ne0,\ \kappa=0,\ K_1=0,\ ({\cal I}_1)$& $\Sigma_{28}$\\[0.7mm]
\cline{3-5} \rule{0pt}{2.8ex}   &     &$ \!{4\choose 1}p\!+\!{0\choose
1}q^c\!+\!{0\choose 1}r^c$
                   & $ \mu_{0,1,2,3}\!=\!0,\ \mu_4\ne0,\ \kappa\ne0,\ ({\cal I}_2)$& $\Sigma_{29}$\\[0.7mm]
\cline{3-5} \rule{0pt}{2.8ex}   &     &$ \!{0\choose 1}p\!+\!{2\choose
1}q^c\!+\!{2\choose 1}r^c$
                   & $ \mu_{0,1,2,3}\!=\!0,\ \mu_4\ne0,\ \kappa=0,\ ({\cal I}_2)$& $\Sigma_{30}$\\[0.7mm]
\cline{2-5} \rule{0pt}{2.8ex} 4 &     &$ {4\choose 1}p+{0\choose 2}q $
                   & $ \mu_{0,1,2,3}\!=\!0,\ \mu_4\ne0,\ \kappa\ne0,\ ({\cal I}_3)$& $\Sigma_{31}$\\[0.7mm]
\cline{3-5} \rule{0pt}{2.8ex}   &     &$ {3\choose 1}p+{1\choose 2}q$
                   & $ \mu_{0,1,2,3}\!=\!0,\ \mu_4\ne0,\ \kappa=L=0,\  \kappa_1\ne0,\  ({\cal I}_3)$& $\Sigma_{32}$\\[0.7mm]
\cline{3-5} \rule{0pt}{2.8ex}   & 2 &$ {2\choose 1}p+{2\choose 2}q$
                   & $ \mu_{0,1,2,3}\!=\!0, \mu_4\ne0, \kappa\!=\!L=\!\kappa_1\!=\!0, K_1=0,  ({\cal I}_3)$& $\Sigma_{33}$\\[0.7mm]
\cline{3-5} \rule{0pt}{2.8ex}   &   &$ {1\choose 1}p+{3\choose 2}q$
            & $ \mu_{0,1,2,3}\!=\!0,\mu_4\ne0, \kappa\!=\!L\!=\!\kappa_1\!=\!0, K_1\ne0,  ({\cal I}_3)$& $\Sigma_{34}$\\[0.7mm]
\cline{3-5} \rule{0pt}{2.8ex}   &    &$ {0\choose 1}p+{4\choose 2}q $
                   & $ \mu_{0,1,2,3}\!=\!0,\ \mu_4\ne0,\ \kappa=0,\ L\ne0,\ ({\cal I}_3)$& $\Sigma_{35}$\\[0.7mm]
\cline{2-5}
\rule{0pt}{2.8ex}  &   3 &$ {4\choose 3}p $   & $ \mu_{0,1,2,3}=0,\ \mu_4\ne0,\ ({\cal I}_4)$& $\Sigma_{36}$\\[0.7mm]
\hline
\end{tabular}
\end{table}

\BProof We need to examine the four distinct cases corresponding to the
canonical forms $(\SSS_I)-(\SSS_{I\!V})$, respectively.
\vspace{-4mm}
\subsection{Systems of type $\SSS_I$}
\vspace{-3mm} For  systems $(\SSS_I)$ we have $\mu_0=gh(g+h-1)$
and for $\mu_0\ne0$ according to Lemma \ref{lm_1} we have
$\Delta_S=0$ and, hence, we obtain a system of the class
$\Sigma_1$ (see Table 2).

Let us consider now $\mu_0=0$. In this case we have $gh(g+h-1)=0$ and without
loss of generality we may assume $g=0$. Indeed, if $h=0$ (respectively,
$g+h-1=0$) we can apply the  linear transformation which will replace the
straight line $y=0$ with $ x=0$ (respectively, $y=0 $ with $y=x$).\ Let $g=0$.
By using the translation
 $x=x_1+(f+eh)/2,$ $ y=y_1+e/2$
we may assume $e=f=0$. In this way the system $(\SSS_I)$ will be brought    to
the following canonical form:
\be
  \dot x=k +cx +dy + (h-1) xy,\quad  \dot y=l -xy +hy^2,
\label{2s4}
\ee
for which we have
$$
\mu_1=ch(1-h)y,\quad \kappa=64h(1-h),\quad  K=2h(h-1)y^2.
$$
For $\mu_1\ne0$, from Lemma~\ref{lm_1}
we obtain $\Delta_S=1$ which leads us to the case $\Sigma_5$.

Considering $\mu_1=0$ we shall examine two cases: $\kappa\ne0$ and $\kappa=0$.
\vspace{-4mm}
\subsubsection{Case $\kappa\ne0$}
\vspace{-3mm}
As the condition $\kappa\ne0$ is equivalent to condition $K\ne0$, according to
Proposition \ref{prop:p2,q2}  we conclude that
 $\Supp D_S(P ,Q ;Z)$ contains exactly one point $p=[1:0:0]$ since
 \mbox{$\gcd(p_2,q_2)=y$}.  By  Lemma \ref{lm_1} its  multiplicity
$I_p(P ,Q )$
 depends of the number of vanishing comitants $\mu_i(a,x,y)$. In this way
we obtain that a quadratic system belongs to the set $\Sigma_{10}$
(respectively $\Sigma_{18}$;\ $\Sigma_{26}$) for $\mu_{0,1}=0,$
$\mu_2\ne0$ (respectively
for $\mu_{0,1,2}=0,\mu_3\ne0$; $\mu_{0,1,2,3}=0$, $\mu_4\ne0$).
\ \ We use the compact notation $\mu_{0,1,2}=0$ for $\mu_0=\mu_1=\mu_2=0$.
\vspace{-5mm}
\subsubsection{Case $\kappa=0$}\label{kapp:0}
\vspace{-4mm}
In this case $h(h-1)=0$ and analogously to the previous case, without loss of
the generality we may assume $h=0$. Thus, for system (\ref{2s4}) we obtain:
\beq
&&\mu_0=\mu_1=0,\qquad \mu_2=-cdxy,\qquad \mu_3=(k-l)(dy-cx)xy,\\
&& \mu_4=-xy[lc^2x^2 -(k-l)^2xy +2lcd xy+ld^2y^2],\quad K_1= -xy(cx+dy).
\eeq
So, if $\mu_2\ne0$ taking into consideration Remark \ref{rm1} and
the value of the comitant $\mu_4$, we obtain the case
$\Sigma_{11}$ in Table 2.

If $\mu_2=0$ and $\mu_3\ne0$ then $cd=0,\ c^2+d^2\ne0$ and clearly we
arrive at the case $\Sigma_{19}$.

Let us now suppose that the conditions $\mu_2=\mu_3=0$ hold.
\vspace{-5mm}
\paragraph{$K_1\ne0$.} Then
$c^2+d^2\ne0$ and from $\mu_3=0$ we obtain $k=l$ which yields
either $\mu_4= -ld^2xy^3$ (for $c=0$) or $\mu_4= -lc^2x^3y$ (for
$d=0$). Both these cases lead us to the case $\Sigma_{27}$ in
Table 2. \vspace{-4mm}
\paragraph{$K_1=0$.} In this case it follows at once that
 $c=d=0$ and, hence, $\mu_4=4(k-l)^2x^2y^2$.
Thus taking into consideration  Remark \ref{rm1} we obtain the case
$\Sigma_{28}$.
\vspace{-4mm}
\subsection{Systems of type $(\SSS_{I\!I})$}
\vspace{-3mm}
For a canonical system $(\SSS_{I\!I})$ we obtain
\beq
&& \mu_0=-h[ g^2+(h+1)^2], \quad \kappa=-64\left[g^2+(h+1)(1-3h)\right],\\
&& K=2(g^2+h+1)x^2 +4ghxy+2h(h+1)y^2
\eeq
and for $\mu_0\ne0$  according to Lemma \ref{lm_1} we have
$\Delta_S=0$. Thus we obtain the case $\Sigma_2$ in Table 2.

Let us consider now $\mu_0=0$, i.e. $h[ g^2+(h+1)^2]=0$.
\vspace{-5mm}
\subsubsection{Case $\kappa\ne0$}
\vspace{-4mm}
 In this case we have $h=0$ and
since  the condition $\kappa\ne0$ is equivalent to the condition $K\ne0$,
according to Proposition  \ref{prop:p2,q2}, $\Supp D_S(P ,Q ;Z)$ contains only
one point,  namely the real one. By Lemma \ref{lm_1} its multiplicity depends
of the number of the vanishing comitants $\mu_i$. Therefore
 the quadratic system belongs to the set $\Sigma_{6}$
(respectively \ $\Sigma_{12}$;\ $\Sigma_{20}$;\ $\Sigma_{29}$)\ for $\mu_1\ne0$
\ (respectively\ for
$\mu_1=0,\mu_2\ne0$;\ $\mu_{1,2}=0,\mu_3\ne0$;\ $\mu_{1,2,3}=0$, $\mu_4\ne0$).
\vspace{-5mm}
\subsubsection{Case $\kappa=0$}
\vspace{-4mm}
The conditions $\mu_0=\kappa=0$ yield $g=0,\ h=-1$ and translating the origin
of coordinates at the point  $(e/4, f/4)$ the system $(\SSS_{I\!I})$  will be
brought to the form
\be
  \dot x=k +cx +dy,\qquad  \dot y=l -x^2 -y^2,
\label{2s5}
\ee
\vspace{-2mm}
for which\\
\vspace{-9mm}
\beq
&&\mu_0=\mu_1=0,\qquad \mu_2=(c^2+d^2)(x^2+y^2),\\
&&\mu_4=  \left ({x}^{2}+{y}^{2}\right )\left [(k^2-c^2l)x^2-2\,cdlxy
  +(k^2-d^2l)y^2\right].
\eeq
 Thus, according to the Remark \ref{rm1}, for $\mu_2\ne0$ we obtain the case $\Sigma_{13}$.

Let us admit that condition $\mu_2=0$ is satisfied. Then $c=d=0$ and for
systems (\ref{2s5}) we have $\mu_3=0,$ $\mu_4=k^2(x^2+y^2)^2$. This leads us to
the case $\Sigma_{30}.$
\vspace{-4mm}
\subsection{Systems of type $(\SSS_{I\!I\!I})$}
\vspace{-3mm}
For canonical systems $(\SSS_{I\!I\!I})$ one can  calculate
\beq
&&\mu_0= gh^2,\quad \kappa= -64h^2,\quad
 K=2\left[g(g-1)x^2+ 2ghxy+ h^2y^2\right].
\eeq
It is quite clear that for $\mu_0\ne0$ we have $\Delta_S=0$ and this
leads us to the case $\Sigma_3$.

Suppose $\mu_0=0$. We examine the two cases: $\kappa\ne0$ and $\kappa =0$.
\vspace{-5mm}
\subsubsection{Case $\kappa\ne0$}
\vspace{-4mm}
Then $h\ne0$ which yields $g=0$ and thus for the
systems $(\SSS_{I\!I\!I})$ we have $\gcd(p_2,q_2)=y$. So, taking
into consideration the Remark \ref{rm1} and the fact that for the
systems $(\SSS_{I\!I\!I})$ the polynomial $C_2(x,y)=x^2y$ we
obtain the case $\Sigma_7$ if $\mu_1\ne0$.

On the other hand the condition $h\ne0$ implies $K\ne0$. Hence, by Proposition
\ref{prop:p2,q2} and Lemma~\ref{lm_1}, $\Supp D_S(P ,Q ;Z)$ contains
 exactly one point $[1:0:0]$ of the  multiplicity $(\Delta_S, 1)$.
 Consequently we conclude that the quadratic system belongs to the set $\Sigma_{14}$\
(respectively,\ $\Sigma_{21}$;\ $\Sigma_{31}$) for
$\mu_1=0,\mu_2\ne0$\ (respectively,\ $\mu_{1,2}=0,\mu_3\ne0$;\
$\mu_{1,2,3}=0$, $\mu_4\ne0$).
\vspace{-5mm}
\subsubsection{Case $\kappa=0$}
\vspace{-4mm}
In this case $h=0$ and for systems
$(\SSS_{I\!I\!I})$  with \quad $p_2=gx^2$,\ $q_2=(g-1)xy$\quad we
have
$$
\mu_0=0,\quad \mu_1= dg(g-1)^2x, \quad L=8gx^2,
$$
and $\gcd(p_2,q_2)=x$. By Lemma~\ref{lm_1} for $\mu_1\ne0$ the quadratic
systems belong to the set $\Sigma_8$.

Supposing $\mu_1=0$ we shall  consider two subcases: $L\ne0$ and $L=0$.
\vspace{-5mm}
\paragraph{Subcase $L\ne0.$}
Then $g\ne0$ and hence $\gcd(p_2,q_2)=x$ for $g\ne1$ and  $\gcd(p_2,q_2)=x^2$
for $g=1$. Hence in both cases by Proposition \ref{prop:p2,q2} and
Lemma~\ref{lm_1}, $\Supp D_S(P ,Q ;Z)$ contains
 exactly one point $[0:1:0]$  whose multiplicity
depends of the number of vanishing comitants $\mu_i(a,x,y)$. Therefore we
conclude that the quadratic systems  belong to the set $\Sigma_{16}$
(respectively $\Sigma_{24}$;\ $\Sigma_{35}$) for $\mu_2\ne0$ (respectively
$\mu_2=0,\mu_3\ne0$; $\mu_{2,3}=0$, $\mu_4\ne0$).
\vspace{-5mm}
\paragraph{Subcase $L=0$.}
For the systems $(\SSS_{I\!I\!I})$ we have $g=0$ and applying the translation
of the phase plane (to obtain $e=f=0$) these systems can be brought to the
form
\be
  \dot x=k +cx+dy,\qquad  \dot y=l - xy.
\label{2s6}
\ee
For the systems (\ref{2s6}) we have $\mu_0=\mu_1=0$ and
$$
 \mu_2=-cdxy,\  \mu_3=-kxy(cx-dy),\  \kappa_1=-32d,\
\mu_4=-xy\left[c^2lx^2+ (2cdl-k^2)xy +d^2ly^2\right].
$$
 So, if $\mu_2\ne0$ by  the Remark \ref{rm1} and Lemma \ref{lm_1}
 the systems (\ref{2s6}) belong to the class $\Sigma_{15}$.

 Let us suppose that the condition $\mu_2=0$ holds.
\vspace{-5mm}
\subparagraph{} If $\kappa_1\ne0$ then $d\ne0$ which implies $c=0$.
Then $\mu_3=dkxy^2$ and taking into consideration the factorization of the
comitant $\mu_4$, we obtain the case $\Sigma_{22}$ for $\mu_3\ne0$ and the case
$\Sigma_{32}$ for $\mu_3=0$, $\mu_4\ne0$.
\vspace{-5mm}
\subparagraph{} Let us suppose $\kappa_1=0$. Then $d=0$
and for the system (\ref{2s6}) we obtain
$$
\mu_3=-ckx^2y,\qquad \mu_4= -x^2y(c^2lx-k^2y),\qquad  K_1=-cx^2y.
$$
Therefore, if $\mu_3\ne0$ by   Remark \ref{rm1} and Lemma \ref{lm_1} the
systems (\ref{2s6}) belong to the class $\Sigma_{23}$. If $\mu_3=0$ we obtain
$ck=0$ and we need to distinguish two cases: $K_1\ne0$ and $K_1=0$.

The condition $K_1\ne0$ yields $c\ne0$ and, hence, $k=0$. This leads us
to the case $\Sigma_{34}$. If $K_1=0$ then $c=0$ and we obtain
the case $\Sigma_{33}$.
\vspace{-4mm}
\subsection{Systems of type $(\SSS_{I\!V})$}
\vspace{-3mm} Note that for systems of the type $(\SSS_{I\!V})$ we
have $D_S(C ,Z)=3q$. So, Supp$D_S(P ,Q ;Z)$ could contain only the
point $[0:1:0]$. By Lemma \ref{lm_1} its multiplicity depends of
the number of the vanishing comitants $\mu_i$. Therefore we obtain
that the quadratic system belongs to the set $\Sigma_{4}$
(respectively\ $\Sigma_9$;\  $\Sigma_{17}$;\ $\Sigma_{25}$;\
$\Sigma_{36}$) for
 $\mu_0\ne0$ (respectively for
$\mu_0=0,\mu_1\ne0$;\ $\mu_{0,1}=0,\mu_2\ne0$;\ $\mu_{0,1,2}=0$, $\mu_3\ne0$;\
$\mu_{0,1,2,3}=0$, $\mu_4\ne0$).

As all cases are examined, Theorem \ref{th_1} is proved. \EProof
\vspace{-4mm}
\section{Divisors encoding the topology of singularities at
infinity}
\vspace{-3mm}
     We  now  need  to  consider  the  topological  types  of  the
singularities  at  infinity of quadratic systems.
For this we shall
introduce  a  third  divisor at infinity:
\bdf\label{df4_1} We call index divisor on the real line at infinity of
$\mathbb R^2$, associated to a real system $(S)$ such that $Z\nmid C $,
  the expression $\sum
i(w)w$ where $w$ is a singular point on the  line  at infinity  $Z  =  0$ of
the system $(S)$ and $i(w)$ is the topological index (cf. \cite{Lib_DS}) of
$w$, i.e. $i(w)$ is the topological index of one of the two opposite singular
points $w,$ $w'$ of ${\cal V}$ on $S^2$.
\edf
\brm\label{rm4_1} This is a well defined divisor which could be
extended trivially to a divisor $\sum j(w)w,$ \ $w\in \{ Z=0\}$ on the line at
infinity $Z=0$ of $\mathbb C^2$ by letting
$$
 j(w)=\left\{\ba{lll} i(w) & if & w\in \PP^2(\mathbb R)\\
                      0    & if & w\in \PP^2(\mathbb C)\setminus \PP^2(\mathbb R),
              \ea\right.
$$
where we identify $\PP^2(\mathbb R)$ with its image via the inclusion
\mbox{$\PP^2(\mathbb R) \hookrightarrow \PP^2(\mathbb C)$} induced by $\mathbb
R \hookrightarrow \mathbb C$.
\erm

\bnot\label{not4_1} We denote  by $I(S)$ the  above divisor on
$Z=0$ in $\PP^2(\C)$, i.e. \mbox{$I(S)=\sum j(w)w$}. \enot

 \bnot\label{not4_2} We denote by  $N_\mathbb C(S)$ (respectively, by $N_\mathbb R(S)$)  the total
 number of
distinct singular points, be they real or complex (respectively, real),
 on the line at infinity $Z  =  0$ of the complex (respectively, real) foliation
 with singularities associated to $(S)$.
\enot

We need to see how the divisor $I(S)=\sum j(w)w$ and the divisors $D_S(P ,Q ;Z)
= \sum I_w(P ,Q )w$ and $D_S(C ,Z) = \sum I_w(C ,Z)p $ constructed in Section 3
are combined. For this we shall fuse these three divisors on the complex line
at infinity into just one but with the values in the  abelian group $\Z^3$:

\bnot Let us consider the following divisor with the value in $\Z^3$ on $Z=0$:
$$
       {\cal D}_S = \sum_w\left(I_w(C ,Z),\, I_w(P ,Q ),\, j(w)  \right)w
$$
where $w$ belongs to the line $Z=0$ of the complex projective plane.
\enot

We cannot detect the multiplicities of the singularities at infinity of a
system $S(\lambda)$ for the parameter value $\lambda$ from just the phase
portrait of $S(\lambda)$.
  On the other hand ${\cal D}_{S(\lambda)}$ has dynamic qualities since it
gives us some information about what could happen to the phase portraits
 in the neighbourhood of $\lambda$. For example if $w\in \{Z=0\}$ and if $I_w(P ,Q )=2$
 for $S(\lambda_0)$,
 then we know that in the neighbourhood of $\lambda_0$ the phase portraits
of $S(\lambda)$ will have 2 finite points arising from $w$ in the neighbourhood
of $w$.

We denote by ${\cal H'}$ and ${\cal H}$ the following sets:
$$
{\cal H'}=\left\{ X^2+Y^2+Z^2=1 \right| \ Z\ge0\left.\right\},\quad {\cal
H}=\left\{ X^2+Y^2+Z^2=1 \right| \ Z>0\left.\right\}.
$$
For $(S)$ in \QS satisfying  the hypothesis of Theorem \ref{th_1} let
$\sigma(S)$ be the set of all $n_\infty=2N_R(S)$ real singular points at
infinity considered on the equator $S^1$ of the Poincar\'e sphere.

We consider the function $n_{sect}:\ \sigma(S)\longrightarrow\ \N$
where $n_{sect}(w)$ is the number of distinct local sectors of the
point $w\in S^1$ on ${\cal H}$.

Let $w\in \sigma(S)$ and let  $\rho(S)=(w_1,w_2,...,w_{n_\infty})$ be the
ordered sequence of singularities of $S$ on $S^1$, enumerated when $S^1$ is
described in the positive sense and such that $w_1=w$.

Let $ O_S(w)=\left(n_{\it sect} (w_1), n_{\it sect}(w_2),...,
n_{\it sect}(w_{n_\infty})\right).  $ Then we have:
$$
O_S(w_i)=\left(n_{\it sect}(w_i),n_{\it sect}(w_{i+1}),...,n_{\it
sect}(w_{n_\infty}), n_{\it sect}(w_1),...,n_{\it
sect}(w_{i-1})\right).
$$
 \bnot\label{not4_3}
We denote by $O(S)$ anyone of the sequences $O_S(w_i)$.
 \enot
 \bnot\label{not4_3a}
 We denote  by $\max(n_{\it sect})$ the maximum value of the function $n_{\it sect}$,
 by  $N_{\max}(n_{\it sect})=\#\{w\in S^1\,| n_{\it sect}(w)=\max(n_{\it
 sect})\}$ and by
  $N_{\it hsect}(S)$ the total number of hyperbolic sectors in ${\cal H}'$ of
singularities at infinity of a  system $(S)\in$ {\bf QS}$_{\bf ess}$.
 \enot

\bdf Let $h_1(w_1)$ and $h_2(w_2)$ be two distinct hyperbolic
sectors of singularities  at infinity $w_1$, $w_2$ of a system
$(S)\in $ {\bf QS}$_{\bf ess}$. (i) We say that $h_1(w_1)$ and
$h_2(w_2)$ are finitely adjacent if $w_1=w_2=w$   and the two
sectors $h_1(w_1)$ and $h_2(w_2)$ have a common border which is a
separatrix of
$w$ in the finite plane.\\
\indent (ii) We say that $h_1(w_1)$ and $h_2(w_2)$ are adjacent at
infinity if $w_1$ and $w_2$ are opposite points of $S^1$ and $w_1$
(also $w_2$) as a point of $S^2$ has two hyperbolic sectors with a
common border, part of the  equator. \edf

 \bnot\label{not4_4} We shall use the following notation
 $$
 N_{\it hsect}^{f\infty a}= ( N_{\it hsect}^{f-a},\, N_{\it
 hsect}^{\infty-a}\,),
 $$
where $N_{\it hsect}^{f-a}$ (respectively $N_{\it
 hsect}^{\infty-a}$)   is the total number of finitely
adjacent couples of hyperbolic sectors (respectively adjacent at infinity).
\enot
\section{Classification
of qua\-dra\-tic differential systems according to their behavior in the
neighborhood of infinity} \label{Sect_class}
 The study of the geometry of the systems yields a simpler set of
algebraic invariants  than those used  in \cite{Nik_Vlp}.  We
refine here the invariants which appeared in \cite{Nik_Vlp} so as
to reveal  the geometry of the systems.

We now need to relate the geometrical invariants defined
in the previous section to their algebraic counterparts,
i.e. the comitants and algebraic invariants.

To do this we construct below the $GL$-comitants which we need, by using the
following basic  ones:
\beq
  && C_i=yp_i(x,y)-xq_i(x,y),\ i=0,1,2, \\
  && D_i=\frac{\partial}{\partial x}p_i(x,y)+
        \frac{\partial}{\partial y}q_i(x,y),\ i=1,2,\quad J_1=Jacob(C_0,D_2),\\
&&  J_2=Jacob(C_0,C_2),\quad J_3=Discrim(C_1),\quad
    J_4=Jacob(C_1,D_2).
\eeq
Using the comitants (\ref{com1}) and (\ref{com2}) we  constructed
in Sections 4 and 5 we  define the following new polynomials:
\be\label{N,R...}
\bal
&  N=K+H,\quad R=  L +8\,K,\quad \kappa_2=-J_1,\quad \xi=M-2K,\\
&  K_2=  4\, Jacob(J_2,\xi)+3\, Jacob(C_1,\xi) D_1-\xi(16J_1+3J_3 +3D_1^2),\\
&  K_3=  2C_2^2(2J_1-3 J_3) + C_2(3C_0K-2C_1J_4) + 2K_1(C_1D_2+3K_1).
\eal
\ee
All these polynomials are  $GL$-comitants, being obtained from simpler
 $GL$-comitants.

In the statement of the next Theorem  {\it Figure} $ j$ \ for
\mbox{$j=${\it 1,...,40}}\ will denote a phase portrait in the
vicinity of infinity of a quadratic system in {\bf QS}$_{\bf
ess}$. The notation for the figures  in \cite{Nik_Vlp} was {\sc
Fig} $ j,$ \ $j=1,...,40$. The correspondence between the two
notations is indicated in columns 6 and 7 in Table 3.

 In our next Theorem we relate the geometry at infinity of quadratic systems
 with algebraic and geometric invariants.

\bth\label{th_2}{\bf[The classification theorem]}
    We consider here the family {\bf QS}$_{\bf ess}$  of all systems $(S)$ in \QS which
are essentially quadratic,  i.e. $\gcd(P,Q)=1$ and $Z\nmid C$.

 {\bf A.} The phase portraits in the vicinity of infinity of the class {\bf QS}$_{\bf ess}$
 are classified topologically by the integer-valued affine invariant
 ${\cal J}= (O,\,N_{\it hsect},\, N_{\it hsect}^{f\infty a}\,)$
which expresses geometrical properties of the systems, e.g. number
of real singularities, number of their sectors and the way in
which these numbers are concatenated, etc. The classification
appears in Table 3 with the corresponding phase portraits  in
Table 5, where they are listed for each value of $N_{\R}(S)$ in
order of increasing topological complexity.

 {\bf B.} The geometrical properties in the neighbourhood of infinity of quadratic systems $(S)$
in {\bf QS}$_{\bf ess}$ are   expressed in terms of algebraic
invariants and comitants as indicated in Table  4, which  contains
the full information regarding multiplicities and indices of the
singularities at infinity for all quadratic differential systems
in {\bf QS}$_{\bf ess}$. The conditions appearing in the last
column of  Table 4 are affinely invariant.
\eth

\begin{table}
\begin{center}
\begin{tabular}{|c|c|c|c|c|c|c|c|}
\multicolumn{8}{r}{\bf Table 3} \\[2mm]
\hline \raisebox{-0.7em}[0pt][0pt]{$N_\R(S)$}
&\raisebox{-0.7em}[0pt][0pt]{$\max(n_{\it sect})$} &
 \raisebox{-0.7em}[0pt][0pt]{$N_{\max}(n_{\it sect})$} &
 \raisebox{-0.7em}[0pt][0pt]{$O(S)$}& \raisebox{-0.7em}[0pt][0pt]{ $N_{\it hsect}$}
 & \multicolumn{2}{c|}{\# of Figures} &\raisebox{-0.7em}[0pt][0pt]{ $ N_{\it hsect}^{f\infty a}$} \\
\cline{6-7}
      &  &       &  &   & New  & Old  &     \\[0.45mm]
 \hline\hline\rule{0pt}{2.65ex}
      & 1 &  6   & (1,1,1,1,1,1) & 0 & 1  & 2  &     \\[0.45mm]
\cline{2-7}\rule{0pt}{2.65ex}
      &   &  1   & (2,1,1,1,1,1) & 2& 2  & 4  &     \\[0.45mm]
\cline{3-7}\rule{0pt}{2.65ex}
      &   &      & (2,2,1,1,1,1) &4 & 3  & 7  &     \\[0.45mm]
\cline{4-7}\rule{0pt}{2.65ex}
  3   & \raisebox{-0.7em}[0pt][0pt]{2}  &  2   & (2,1,2,1,1,1)& 4 & 4  & 6  &     \\[0.45mm]
\cline{4-7}\rule{0pt}{2.65ex}
      &   &      & (2,1,1,2,1,1) & 4 & 5  & 1  &     \\[0.45mm]
\cline{3-7}\rule{0pt}{2.65ex}
      &   &  3   & (2,2,1,1,2,1) & 6 & 6  & 5  &     \\[0.45mm]
\cline{3-7}\rule{0pt}{2.65ex}
      &   &  4   & (2,2,1,2,2,1) & 8  & 7  & 3  &     \\[0.45mm]
\cline{1-7}\rule{0pt}{2.65ex}
       &   &      &              & 2 &  8 & 22  &     \\[0.45mm]
\cline{5-7}\rule{0pt}{2.65ex}
       & 1 &  4  & (1,1,1,1)     & 1 & 9 & 12  &     \\[0.45mm]
\cline{5-7}\rule{0pt}{2.65ex}
       &  &      &               & 0 & 10  & 18  &     \\[0.45mm]
\cline{2-7}\rule{0pt}{2.65ex}
       &   &      &              & 3 &  11  & 15  &     \\[0.45mm]
\cline{5-8}\rule{0pt}{2.65ex}
       &   &  \raisebox{-0.7em}[0pt][0pt]{1}   & \raisebox{-0.7em}[0pt][0pt]{(2,1,1,1)}    & \raisebox{-0.7em}[0pt][0pt]{2} &  12  & 26  & (2,0)    \\[0.45mm]
\cline{6-8}\rule{0pt}{2.65ex}
       &   &     &               &    & 13  & 16  &  (0,2)   \\[0.45mm]
\cline{5-8}\rule{0pt}{2.65ex}
       &   &     &               & 1  & 14   & 23  &     \\[0.45mm]
\cline{3-7}\rule{0pt}{2.65ex}
       &   &    & (2,2,1,1)    &  3  &  15  & 29  &     \\[0.45mm]
\cline{4-7}\rule{0pt}{2.65ex}
       &\raisebox{-0.7em}[0pt][0pt]{ 2}   &      &  & 5 &  16  & 13  &    \\[0.45mm]
\cline{5-7}\rule{0pt}{2.65ex}
       &  &  2   & \raisebox{-0.7em}[0pt][0pt]{(2,1,2,1)}            &4  & 17  & 20  &    \\[0.45mm]
\cline{5-8}\rule{0pt}{2.65ex}
    \raisebox{-0.7em}[0pt][0pt]{ 2}   &   &      &            &\raisebox{-0.7em}[0pt][0pt]{2}  & 18  & 8  & (0,1)   \\[0.45mm]
\cline{6-8}\rule{0pt}{2.65ex}
       &     &      &           &   & 19   & 21  & (0,0)    \\[0.45mm]
\cline{3-8}\rule{0pt}{2.65ex}
       &   &  \raisebox{-0.7em}[0pt][0pt]{3}   &
       \raisebox{-0.7em}[0pt][0pt]{(2,1,2,2)}& \raisebox{-0.7em}[0pt][0pt]{ 4} & 20   & 10  & (2,2)     \\[0.45mm]
\cline{6-8}\rule{0pt}{2.65ex}
       &   &      &            &   &  21  & 25  & (2,0)   \\[0.45mm]
\cline{3-8}\rule{0pt}{2.65ex}
      &   &  4   & (2,2,2,2) & 6 &  22  & 9  &    \\[0.45mm]
\cline{2-7}\rule{0pt}{2.65ex}
       &   &      & \raisebox{-0.7em}[0pt][0pt]{(3,1,1,1)}   & 4 &  23  & 11  &      \\[0.45mm]
\cline{5-7}\rule{0pt}{2.65ex}
       &   & \raisebox{-0.7em}[0pt][0pt]{1}   &     & 3 & 24   & 28  &      \\[0.45mm]
\cline{4-7}\rule{0pt}{2.65ex}
       &  &       &  {(3,1,2,1)}&  4 &  25  & 24  &     \\[0.45mm]
\cline{4-7}\rule{0pt}{2.65ex}
       & 3  &       &  {(3,2,1,2)}& 5  & 26   & 14  &     \\[0.45mm]
\cline{3-7}\rule{0pt}{2.65ex}
       &  &      &  \raisebox{-0.7em}[0pt][0pt]{(3,1,3,1)}&  6 & 27   & 19  &     \\[0.45mm]
\cline{5-7}\rule{0pt}{2.65ex}
       &  &   2    &   & 2 & 28   & 27  &     \\[0.45mm]
\cline{4-7}\rule{0pt}{2.65ex}
       &  &       &  {(3,2,3,2)}&  6  &  29  & 17  &     \\[0.45mm]
\hline

\end{tabular}
\end{center}
\end{table}

\begin{table}
\quad\begin{tabular}{|c|c|c|c|c|c|c|c|}
\multicolumn{8}{r}{\bf Table 3}({\it continued}) \\[2mm]
\hline \raisebox{-0.7em}[0pt][0pt]{$N_\R(S)$}
&\raisebox{-0.7em}[0pt][0pt]{$\max(n_{\it sect})$} &
 \raisebox{-0.7em}[0pt][0pt]{$N_{\max}(n_{\it sect})$} &
 \raisebox{-0.7em}[0pt][0pt]{$O(S)$}& \raisebox{-0.7em}[0pt][0pt]{ $N_{\it hsect}$}
 & \multicolumn{2}{c|}{\# of Figures} &\raisebox{-0.7em}[0pt][0pt]{ $ N_{\it hsect}^{f\infty a}$} \\
\cline{6-7}
      &  &       &  &   & New  & Old  &     \\[0.30mm]
 \hline\hline\rule{0pt}{2.60ex}
     & 1 &  2   & (1,1) & 0  & 30  & 30  &     \\[0.30mm]
\cline{2-7}\rule{0pt}{2.60ex}
      &  &      &       & 2 & 31  & 32  &     \\[0.30mm]
\cline{5-7}\rule{0pt}{2.60ex}
      &   &  1   & (2,1) & 1 & 32 & 34  &     \\[0.30mm]
\cline{5-7}\rule{0pt}{2.60ex}
       & \raisebox{-0.7em}[0pt][0pt]{2}  &    &     & 0 & 33   & 38  &     \\[0.30mm]
\cline{3-7}\rule{0pt}{2.60ex}
       &   &     &   & 4 &  34  & 31  &     \\[0.30mm]
\cline{5-8}\rule{0pt}{2.60ex}
   1    &   & 2  &  (2,2)   & \raisebox{-0.7em}[0pt][0pt]{2} &  35  & 40  & (2,0)    \\[0.30mm]
\cline{6-8}\rule{0pt}{2.60ex}
       &   &    &     &   & 36   & 39  & (0,2)    \\[0.30mm]
\cline{2-8}\rule{0pt}{2.60ex}
       & &     &  \raisebox{-0.7em}[0pt][0pt]{(3,1)}  & \raisebox{-0.7em}[0pt][0pt]{2} & 37 & 33  & (2,0)    \\[0.30mm]
\cline{6-8}\rule{0pt}{2.60ex}
       & \raisebox{-0.7em}[0pt][0pt]{3}  &  1   &         &  &  38  & 37  &  (0,0)   \\[0.30mm]
\cline{4-8}\rule{0pt}{2.60ex}
       &   &     &  (3,2)  & 3&  39  & 36  &     \\[0.30mm]
\cline{3-7}\rule{0pt}{2.60ex}
       &   &  2   &  (3,3)  & 4&  40  & 35  &     \\[0.30mm]
\hline
\end{tabular}

\vspace{4mm}
\noindent\begin{tabular}{|l|c|l|}
\multicolumn{3}{r}{\bf Table 4}\\[1.5mm]
\hline \!\! Figures\!\! & Value of ${\cal D}_S$  & Necessary
and sufficient conditions \\
 \hline\hline\rule{0pt}{2.8ex}
 \raisebox{-0.7em}[0pt][0pt]{Fig.\,1}
  & $ (1,0,1)p+(1,0,1)q+(1,0,1)r$ &  $\eta>0,\,  \mu_0<0,\, \kappa>0$\\[0.5mm]
\cline{2-3}\rule{0pt}{2.8ex}
   & $ (1,2,1)p+(1,0,1)q+(1,0,1)r$ &  $\eta>0,\, \mu_{0,1}=0,\, \mu_2<0,\, \kappa>0$\\[0.5mm]
\hline \rule{0pt}{2.8ex}
  & $ (1,1,0)p+(1,0,1)q+(1,0,1)r$ &  $\eta>0,\,  \mu_0=0,\,\mu_1\ne0,\, \kappa>0$\\[0.5mm]
\cline{2-3}\rule{0pt}{2.8ex} {Fig.\,2}
   & $ (1,3,1)p+(1,0,1)q+(1,0,1)r$ &  $\eta>0,\, \mu_{0,1,2}=0,\, \mu_3\ne0,\, \kappa>0$\\[0.5mm]
\cline{2-3}\rule{0pt}{2.8ex}
   & $ (1,2,1)p+(1,1,1)q+(1,0,1)r$ &  $\eta>0,\, \mu_{0,1,2}=\kappa=0,\, \mu_3 K_1<0 $\\[0.5mm]
\hline \rule{0pt}{2.8ex}
 \raisebox{-0.7em}[0pt][0pt]{Fig.\,3}
   & $ (1,1,0)p+(1,1,0)q+(1,0,1)r$ &  $\eta>0,\, \mu_{0,1}=\kappa=0,\, \mu_2L<0 $\\[0.5mm]
\cline{2-3}\rule{0pt}{2.8ex}
   & $ (1,3,0)p+(1,1,0)q+(1,0,1)r$ &  $\eta>0,\, \mu_{0,1,2,3}=\kappa=0,\, \mu_4L<0,\, K_1\ne0 $\\[0.5mm]
\hline \rule{0pt}{2.8ex}
 \raisebox{-0.7em}[0pt][0pt]{Fig.\,4}
   & $ (1,1,0)p+(1,1,0)q+(1,0,1)r$ &  $\eta>0,\, \mu_{0,1}=\kappa=0,\, \mu_2L>0 $\\[0.5mm]
\cline{2-3}\rule{0pt}{2.8ex}
   & $ (1,3,0)p+(1,1,0)q+(1,0,1)r$ &  $\eta>0,\, \mu_{0,1,2,3}=\kappa=0,\, \mu_4L>0,\, K_1\ne0 $\\[0.5mm]
\hline \rule{0pt}{2.8ex}
  & $ (1,0,1)p+(1,0,1)q+(1,0,-1)r$ &  $\eta>0,\,  \mu_0>0$\\[0.5mm]
\cline{2-3}\rule{0pt}{2.8ex}
   & $ (1,2,1)p+(1,0,1)q+(1,0,-1)r$ &  $\eta>0,\, \mu_{0,1}=0,\, \mu_2>0,\, \kappa<0$\\[0.5mm]
\cline{2-3}\rule{0pt}{2.8ex}
 \raisebox{-0.7em}[0pt][0pt]{Fig.\,5}
   & $ (1,4,1)p+(1,0,1)q+(1,0,-1)r$ &  $\eta>0,\, \mu_{0,1,2,3}=0,\, \mu_4\ne0,\, \kappa<0$\\[0.5mm]
\cline{2-3}\rule{0pt}{2.8ex}
   & $ (1,0,1)p+(1,0,1)q+(1,2,-1)r$ &  $\eta>0,\, \mu_{0,1}=0,\, \mu_2>0,\,
   \kappa>0$\\[0.5mm]
\cline{2-3}\rule{0pt}{2.8ex}
   & $ (1,0,1)p+(1,0,1)q+(1,4,-1)r$ &  $\eta>0,\, \mu_{0,1,2,3}=0,\, \mu_4\ne0,\, \kappa>0$\\[0.5mm]
\cline{2-3}\rule{0pt}{2.8ex}
   & $ (1,2,1)p+(1,0,1)q+(1,2,-1)r$ &  $\eta>0,\, \mu_{0,1,2,3}=0,\, \mu_4\ne0,\, \kappa= K_1=0$\\[0.5mm]
 \hline
\end{tabular}
\end{table}

\begin{table}
\noindent\begin{tabular}{|l|c|l|}
\multicolumn{3}{r}{\bf Table 4}({\it continued})\\[2mm]
\hline \!\! Figures\!\! & Value of ${\cal D}_S$  & Necessary
and sufficient conditions \\
 \hline\hline\rule{0pt}{2.8ex}
  & $ (1,1,0)p+(1,0,1)q+(1,0,-1)r$ &  $\eta>0,\,  \mu_0=0,\,\mu_1\ne0,\, \kappa<0$\\[0.5mm]
\cline{2-3}\rule{0pt}{2.8ex} {Fig.\,6}
   & $ (1,3,1)p+(1,0,1)q+(1,0,-1)r$ &  $\eta>0,\, \mu_{0,1,2}=0,\, \mu_3\ne0,\, \kappa<0$\\[0.5mm]
\cline{2-3}\rule{0pt}{2.8ex}
   & $ (1,2,1)p+(1,1,1)q+(1,0,-1)r$ &  $\eta>0,\, \mu_{0,1,2}=\kappa=0,\, \mu_3 K_1>0 $\\[0.5mm]
\hline \rule{0pt}{2.8ex}
 \raisebox{-0.7em}[0pt][0pt]{Fig.\,7}
  & $\! (1,0,\!-1)p\!+\!(1,0,1)q\!+\!(1,0,-1)r\!$ &  $\eta>0,\,  \mu_0<0,\, \kappa<0$\\[0.5mm]
\cline{2-3}\rule{0pt}{2.8ex}
   & $\! (1,2,\!-1)p\!+\!(1,0,1)q\!+\!(1,0,-1)r\!$ &  $\eta>0,\, \mu_{0,1}=0,\, \mu_2<0,\, \kappa<0$\\[0.5mm]
\hline\rule{0pt}{2.8ex}
   & $ (2,2,0)\,p+(1,0,1)\,q $  &  $\!\!\! \ba{l}\eta=0,\, M\ne0,\, \mu_{0,1}=\kappa=\kappa_1=0,\\
                       \hspace{12mm}\mu_2>0,\,L>0,\, K_2<0\ea \!$\\[0.50mm]
\cline{2-3}\rule{0pt}{2.4ex}
 Fig.\,8
   & $\!(2,4,0)\,p+(1,0,1)\,q\!$  &  $\!\!\! \ba{l}\eta=0,\, M\ne0,\, \mu_{0,1,2,3}=\kappa=\kappa_1=0,\\
                          \hspace{12mm}\mu_4>0,\,L>0,\,K=0,\,K_2<0\ea $\\[0.30mm]
\cline{2-3}\rule{0pt}{2.8ex}
   & $\!(2,2,0)\,p+(1,2,1)\,q\!$  &  $\!\!\! \ba{l}\!\eta=0, M\ne0,\, \mu_{0,1,2,3}=\kappa=\kappa_1=0\\
                              \hspace{12mm}\mu_4\ne0,\, L=K_1=0,\,  \kappa_2<0\ea $\\[0.50mm]
\hline\rule{0pt}{2.8ex}
   & $ (2,1,1)\,p+(1,0,1)\,q $  &  $\eta=0, M\mu_1\!\ne\!0,\, \mu_{0}\!=\!\kappa\!=\!0,L\!>\!0,K\!<\!0 $\\[0.50mm]
\cline{2-3}\rule{0pt}{2.8ex} Fig.\,9
   & $ (2,3,1)\,p+(1,0,1)\,q $  &  $\eta=0, M\kappa_1L\!\ne\!0, \mu_{0,1,2}\!=\!\kappa\!=\!0,\mu_3K_1\!<\!0 $\\[0.50mm]
\cline{2-3}\rule{0pt}{2.5ex}
   & $ (2,1,1)\,p+(1,2,1)\,q $  &  $\!\!\!\ba{l}\eta=0,\, M\ne0,\, \mu_{0,1,2}=\kappa=L=0,\\[-1mm]
                              \hspace{12mm}\kappa_1\ne0,\ \mu_3K_1<0\ea $\\[0.40mm]
\hline\rule{0pt}{2.8ex}
 Fig.\,10
   & $ (2,2,2)\,p+(1,0,1)\,q $  &  $\!\!\!\ba{l}\eta=0,\, M\ne0,\, \mu_{0,1}\!=\!\kappa\!=\!\kappa_1\!=\!0,\\
                               \hspace{12mm}\mu_2<0,\,L>0,\, K<0 \ea\!$\\[0.50mm]
\hline\rule{0pt}{2.8ex}
 \raisebox{-0.7em}[0pt][0pt]{Fig.\,11}
   & $ (2,1,1)\,p+(1,1,0)\,q $  &  $\eta=0,\, M\ne0,\, \mu_{0,1}=\kappa=L=0,\,\mu_2\ne0 $\\[0.50mm]
\cline{2-3}\rule{0pt}{2.8ex}
   & $ (2,1,1)\,p+(1,3,0)\,q $  &  $\eta=0, M\ne0, \mu_{0,1,2,3}\!=\!\kappa\!=L\!=\!0,\mu_4\kappa_1\!\ne\!0  $\\[0.50mm]
\hline\rule{0pt}{2.8ex}
 Fig.\,12
   & $\! (2,2,2)\,p+(1,1,0)\,q\! $  &  $\!\!\!\ba{l}\eta=0, M\ne0, \mu_{0,1,2}=\kappa=\kappa_1=0,\\
                         \hspace{12mm} L=0,\,\mu_3 K_1<0 \ea$\\[0.50mm]
\hline\rule{0pt}{2.8ex}
 \raisebox{-0.7em}[0pt][0pt]{Fig.\,13}
   & $ (2,2,1)\,p+(1,0,1)\,q $  &  $\eta=0, M\ne0, \mu_{0,1}=\kappa=0,\mu_2\!\ne\!0,\kappa_1L\!\ne\!0 $\\[0.50mm]
\cline{2-3}\rule{0pt}{2.8ex}
   & $ (2,4,1)\,p+(1,3,0)\,q $  &  $\eta=0, M\mu_4\!\ne\!0, \mu_{0,1,2,3}\!=\!\kappa\!=\!0,\kappa_1L\!\ne\!0  $\\[0.50mm]
\hline\rule{0pt}{2.8ex}
 Fig.\,14
   & $ (2,3,1)\,p+(1,0,1)\,q $  &  $\!\!\!\ba{l}\eta=0,\,M\ne0,\mu_{0,1,2}\!=\!\kappa\!=\!\kappa_1\!=\!0, \\
                                 \hspace{12mm}\mu_3\ne0,\,L>0,\, K<0 \ea\!$\\[0.50mm]
\hline\rule{0pt}{2.8ex}
 Fig.\,15
   & $\!(2,3,1)\,p+(1,1,0)\,q\!$  &  $\eta=0,\, M\mu_4K_1\!\ne\!0,\, \mu_{0,1,2,3}\!=\!\kappa\!=\!\kappa_1\!=\!L\!=\!0\!\!\! $\\[0.50mm]
\hline\rule{0pt}{2.8ex}
 \raisebox{-0.7em}[0pt][0pt]{Fig.\,16}
   & $ \!(2,1,1)\,p+(1,0,-1)\,q \!$  &  $\eta=0, M\mu_1\ne0, \mu_{0}\!=\!\kappa\!=\!0,L\!<\!0,N\!\le\!0 $\\[0.50mm]
\cline{2-3}\rule{0pt}{2.8ex}
   & $\! (2,1,1)\,p+(1,2,-1)\,q \!$  &  $\!\!\!\ba{l} \eta=0,\, M\ne0,\, \mu_{0,1,2}\!=\kappa=L=0,\\
                                \hspace{12mm}\kappa_1\ne0,\,\mu_3K_1>0 \ea $\\[0.50mm]
 \hline
\end{tabular}
\end{table}

\begin{table}
\noindent\begin{tabular}{|l|c|l|}
\multicolumn{3}{r}{\bf Table 4} ({\it continued}\,)\\[1mm]
\hline \!\! Figures\!\! & Value of ${\cal D}_S$  &\hfil Necessary
and sufficient conditions\hfill \\
 \hline\hline\rule{0pt}{2.8ex}
   & $\!(2,2,2)\,p+(1,0,-1)\,q\!$  &  $\eta=0,\, M\ne0,\, \mu_{0,1}\!=\!\kappa\!=\!\kappa_1\!=\!0,\,\mu_2>0,\,L<0 $\\[0.50mm]
\cline{2-3}\rule{0pt}{2.5ex}
 Fig.\,17
   & $\!(2,4,2)\,p+(1,0,-1)\,q\!$  &  $\eta=0,\, M\ne0,\, \mu_{0,1,2,3}\!=\!\kappa\!=\!\kappa_1\!=\!0,\,\mu_4>0,\,L<0 $\\[0.30mm]
\cline{2-3}\rule{0pt}{2.8ex}
   & $\!(2,2,2)\,p+(1,2,-1)\,q\!$  &  $\!\eta\!=\!0, M\!\ne\!0,\, \mu_{0,1,2,3}\!=\!\kappa\!=\!\kappa_1\!=\!L\!=\!K_1\!=\!0, \mu_4\!\ne\!0, \kappa_2\!>\!0\! $\\[0.50mm]
\hline\rule{0pt}{2.8ex}
   & $ (2,0,0)\,p+(1,0,1)\,q $  &  $\eta=0,\, M\ne0,\, \mu_{0}>0 $\\[0.5mm]
\cline{2-3}\rule{0pt}{2.8ex}
   & $ (2,0,0)\,p+(1,2,1)\,q $  &  $\eta=0,\, M\ne0,\, \mu_{0,1}=0,\, \mu_{2}>0,\,\kappa\ne0 $\\[0.5mm]
\cline{2-3}\rule{0pt}{2.8ex}
 \raisebox{-1.9em}[0pt][0pt]{Fig.\,18}
    & $ (2,0,0)\,p+(1,4,1)\,q $  &  $\eta=0,\, M\ne0,\, \mu_{0,1,2,3}=0,\, \mu_{4}\ne0,\,\kappa\ne0 $\\[0.5mm]
\cline{2-3}\rule{0pt}{2.8ex} &
\raisebox{-1.3em}[0pt][0pt]{$(2,4,0)\,p+(1,0,1)\,q $}  &
         $\!\!\! \ba{l} \eta=0,\, M\ne0,\,  \mu_{0,1,2,3}=\kappa=\kappa_1=0,\\[-0.8mm]
         \hphantom{\eta=0,\ }\mu_{4}>0,\, L>0,\, K\ne0,\, R\ge0 \ea  $\\[0.5mm]
\cline{3-3}\rule{0pt}{2.8ex}
 & & $\!\!\! \ba{l} \eta=0,\, M\ne0,\,  \mu_{0,1,2,3}=\kappa=\kappa_1=0,\\[-0.5mm]
         \hphantom{\eta=0,\ }\mu_{4}\!>\!0,\, L>0,\, K\!=\!0,\,K_2\ge0 \ea  $\\[0.5mm]
\hline\rule{0pt}{2.8ex}
 Fig.\,19
   & $ (2,2,0)\,p+(1,0,1)\,q $  &  $\eta=0,\, M\ne0,\, \mu_{0,1}\!=\!\kappa\!=\!\kappa_1\!=\!0,\,\mu_2>0,\,L\!>\!0,\, K_2\!\ge\!0 \!$\\[0.50mm]
\hline \rule{0pt}{2.8ex}
 \raisebox{-0.7em}[0pt][0pt]{Fig.\,20}
   & $ (2,0,0)\,p+(1,1,0)\,q $  &  $\eta=0,\, M\ne0,\, \mu_{0}=0,\,\mu_1\ne0,\,\kappa\ne0 $\\[0.5mm]
\cline{2-3}\rule{0pt}{2.8ex}
   & $ (2,0,0)\,p+(1,3,0)\,q $  &  $\eta=0,\, M\ne0,\, \mu_{0,1,2}=0,\,\mu_3\ne0,\,\kappa\ne0 $\\[0.5mm]
\hline\rule{0pt}{2.8ex}
 Fig.\,21
   & $\! (2,2,0)\,p+(1,1,0)\,q\! $  &  $\eta=0,\, M\ne0,\, \mu_{0,1,2}\!=\!\kappa\!=\!\kappa_1\!=\!L\!=\!0,\,\mu_3 K_1\!>\!0 \!$\\[0.50mm]
\hline\rule{0pt}{2.8ex}
 \raisebox{-0.7em}[0pt][0pt]{Fig.\,22}
   & $ (2,0,0)\,p+(1,0,1)\,q $  &  $\eta=0,\, M\ne0,\, \mu_{0}<0 $\\[0.5mm]
\cline{2-3}\rule{0pt}{2.8ex}
   & $ (2,0,0)\,p+(1,2,1)\,q $  &  $\eta=0,\, M\ne0,\, \mu_{0,1}=0,\, \mu_{2}<0,\,\kappa\ne0 $\\[0.5mm]
\hline\rule{0pt}{2.8ex}
 \raisebox{-0.7em}[0pt][0pt]{Fig.\,23}
   & $\! (2,1,-1)\,p+(1,0,1)\,q \!$  &  $\eta=0, M\!\ne\!0, \mu_{0}\!=\!\kappa\!=\!0,\mu_1\!\ne\!0,L\!>\!0,K\!>\!0\! $\\[0.50mm]
\cline{2-3}\rule{0pt}{2.8ex}
   & $\! (2,3,-1)\,p+(1,0,1)\,q\! $  &  $\eta=0, M\!\ne\!0, \mu_{0,1,2}\!=\!\kappa\!=\!0,\kappa_1L\!\ne\!0,\mu_3K_1\!>\!0\! $\\[0.50mm]
\hline\rule{0pt}{2.8ex}
 Fig.\,24
   & $\!(2,4,0)\,p+(1,0,1)\,q\!$  &  $\eta=0,\, ML\ne0,\, \mu_{0,1,2,3}\!=\!\kappa\!=\!\kappa_1\!=\!0,\mu_4\!<\!0\! $\\[0.50mm]
\hline\rule{0pt}{2.8ex}
 Fig.\,25
   & $\! (2,3,-1)\,p+(1,0,1)\,q\! $  &  $\!\!\!\ba{l} \eta=0,\, M\ne0,\, \mu_{0,1,2}\!=\!\kappa\!=\!\kappa_1\!=\!0,\\
                      \hspace{12mm}\mu_3\ne0,\,L>0,\, K>0 \ea\!$\\[0.50mm]
\hline\rule{0pt}{2.8ex}
 \raisebox{-0.7em}[0pt][0pt]{Fig.\,26}
   & $ \!(2,1,1)\,p+(1,0,-1)\,q \!$  &  $\eta=0, M\!\ne\!0, \mu_{0}\!=\!\kappa\!=\!0,\mu_1\!\ne\!0,L\!<\!0,N\!>\!0 $\\[0.50mm]
\cline{2-3}\rule{0pt}{2.8ex}
   & $\! (2,3,1)\,p+(1,0,-1)\,q \!$  &  $\eta=0, M\!\ne\!0, \mu_{0,1,2}\!=\!\kappa\!=\!\kappa_1\!=\!0,\mu_3\!\ne\!0,L\!<\!0  $\\[0.50mm]
\hline\rule{0pt}{2.8ex}
 Fig.\,27
   & $\! (2,2,-2)\,p+(1,0,1)\,q \!$  &  $\!\!\!\ba{l}\eta=0,\, M\ne0,\, \mu_{0,1}\!=\!\kappa\!=\!\kappa_1\!=\!0,\\
                               \hspace{12mm}\mu_2<0,\,L>0,\, K>0 \ea\!$\\[0.50mm]
\hline\rule{0pt}{2.8ex}
 Fig.\,28
   & $\!(2,4,0)\,p+(1,0,1)\,q\!$  &  $\!\!\!\ba{l}\eta=0,\, M\ne0,\,
                                        \mu_{0,1,2,3}=\kappa=\kappa_1=0,\\[-0.5mm]
                                        \hphantom{\eta=0,\ } \mu_4>0,\,L>0,\, K\ne0,\,R<0 \ea$\\[0.50mm]
\hline\rule{0pt}{2.8ex}
 Fig.\,29
   & $\!(2,2,0)\,p+(1,0,-1)\,q\! $  &  $\eta=0, M\ne0, \mu_{0,1}\!=\!\kappa\!=\!\kappa_1\!=\!0,\mu_2\!<\!0,L\!<\!0 $\\[0.50mm]
 \hline
\end{tabular}
\end{table}

\begin{table}
\noindent\begin{tabular}{|l|c|l|}
\multicolumn{3}{r}{\bf Table 4} ({\it continued}\,)\\[1mm]
\hline \!\! Figures\!\! & Value of ${\cal D}_S$  &\hfil Necessary
and sufficient conditions\hfill \\
 \hline\hline\rule{0pt}{2.8ex}
   & $\!(1,0,1)p+(1,0,0)q^c+(1,0,0)r^c\!$  &  $\eta<0,\,  \mu_{0}>0\! $\\[0.50mm]
\cline{2-3}\rule{0pt}{2.8ex}
   & $\!(1,2,1)p+(1,0,0)q^c+(1,0,0)r^c\!$  &  $\eta<0,\,  \mu_{0,1}=0,\,\mu_2>0,\,\kappa\ne0\! $\\[0.50mm]
\cline{2-3}\rule{0pt}{2.8ex}
   & $\!(1,4,1)p+(1,0,0)q^c+(1,0,0)r^c\!$  &  $\eta<0,\,  \mu_{0,1,2,3}=0,\,\mu_4\ne0,\,\kappa\ne0\! $\\[0.50mm]
\cline{2-3}\rule{0pt}{2.8ex}
   & $\!(1,0,1)p+(1,1,0)q^c+(1,1,0)r^c\!$  &  $\eta<0,\,  \mu_{0,1}=\kappa=0,\,\mu_2\ne0\! $\\[0.50mm]
\cline{2-3}\rule{0pt}{2.8ex}
 Fig.\,30
   & $\!(1,0,1)p+(1,2,0)q^c+(1,2,0)r^c\!$  &  $\eta<0,\,  \mu_{0,1,2,3}=\kappa=0,\,\mu_4\ne0  $\\[0.50mm]
\cline{2-3}\rule{0pt}{2.8ex}
   & $ (3,0,1)\,p$     &  $M=0,\,  \mu_{0}>0 $\\[0.50mm]
\cline{2-3}\rule{0pt}{2.8ex}
   & \raisebox{-0.7em}[0pt][0pt]{$ (3,2,1)\,p$}     &  $M=0,\,  \mu_{0,1}=0,\,\mu_2>0,\, K\ne0,\,K_2<0 $\\[0.50mm]
\cline{3-3}\rule{0pt}{2.8ex}
   &      &  $M=0,\,  \mu_{0,1}=0,\,\mu_2>0,\, K=0 $\\[0.50mm]
\cline{2-3}\rule{0pt}{2.8ex}
   & $ (3,4,1)\,p$     &  $M=0,\,  \mu_{0,1,2,3}=0,\,\mu_4>0,\, K_3\ge0 $\\[0.50mm]
\hline\rule{0pt}{2.8ex}
   & $\!(1,1,0)p+(1,0,0)q^c+(1,0,0)r^c\!$  &  $\eta<0,\,  \mu_{0}=0,\,\mu_1\ne0\! $\\[0.50mm]
\cline{2-3}\rule{0pt}{2.8ex}
 Fig.\,31
   & $\!(1,3,0)p+(1,0,0)q^c+(1,0,0)r^c\!$  &  $\eta<0,\,  \mu_{0,1,2}=0,\,\mu_3\ne0\! $\\[0.50mm]
\cline{2-3}\rule{0pt}{2.8ex}
   & $ (3,3,0)\,p$     &  $M=0,\,  \mu_{0,1,2}=K=0,\,\mu_3K_1\!>\!0,\, K_3\!\ge\!0\! $\\[0.50mm]
 \hline\rule{0pt}{2.8ex}
  \raisebox{-0.7em}[0pt][0pt]{Fig.\,32}
   & $ (3,2,1)\,p$     &  $M=0,\,  \mu_{0,1}=0,\,\mu_2>0,\,K\ne0,\,K_2\ge0 $\\[0.50mm]
\cline{2-3}\rule{0pt}{2.8ex}
   & $ (3,4,1)\,p$     &  $M=0,\,  \mu_{0,1,2,3}=K=0,\,\mu_4>0,\, K_3<0\! $\\[0.50mm]
\hline\rule{0pt}{2.8ex}
 Fig.\,33
   & $ (3,3,2)\,p$     &  $M=0,\,  \mu_{0,1,2}=K=0,\,\mu_3K_1<0\! $\\[0.50mm]
\hline\rule{0pt}{2.8ex}
   & $\!(1,0,-1)p+(1,0,0)q^c+(1,0,0)r^c\!$  &  $\eta<0,\,  \mu_{0}<0\! $\\[0.50mm]
\cline{2-3}\rule{0pt}{2.8ex}
 Fig.\,34
   & $\!(1,2,-1)p+(1,0,0)q^c+(1,0,0)r^c\!$  &  $\eta<0,\,  \mu_{0,1}=0,\,\mu_2<0,\,\kappa\ne0\! $\\[0.50mm]
\cline{2-3}\rule{0pt}{2.8ex}
   & $ (3,0,-1)\,p$     &  $M=0,\,  \mu_{0}<0 $\\[0.50mm]
\hline\rule{0pt}{2.8ex}
 Fig.\,35
   & $ (3,4,1)\,p$     &  $M=0,\,  \mu_{0,1,2,3}=0,\,\mu_4<0\! $\\
\hline\rule{0pt}{2.8ex}
 Fig.\,36
   & $ (3,4,1)\,p$     &  $M=0,\,  \mu_{0,1,2,3}=0,\,\mu_4>0, K\!\ne\!0, K_3\!<\!0\! $\\[0.50mm]
\hline\rule{0pt}{2.8ex}
 \raisebox{-0.7em}[0pt][0pt]{Fig.\,37}
   & $ (3,1,0)\,p$     &  $M=0,\,  \mu_{0}=0,\,\mu_1\ne0 $\\[0.50mm]
\cline{2-3}\rule{0pt}{2.8ex}
   & $ (3,3,0)\,p$     &  $M=0,\,  \mu_{0,1,2}=0,\,\mu_3K\ne0,\, K_3>0\! $\\[0.50mm]
\hline\rule{0pt}{2.8ex}
 Fig.\,38
   & $ (3,3,0)\,p$     &  $M=0,\,  \mu_{0,1,2}=K=0,\,\mu_3K_1\!>\!0,\, K_3\!<\!0\! $\\[0.50mm]
\hline\rule{0pt}{2.8ex}
 Fig.\,39
   & $ (3,3,0)\,p$     &  $M=0,\,  \mu_{0,1,2}=0,\,\mu_3K\ne0,\, K_3<0\! $\\[0.50mm]
\hline\rule{0pt}{2.8ex}
 Fig.\,40
   & $ (3,2,-1)\,p$     &  $M=0,\,  \mu_{0,1}=0,\,\mu_2<0\! $\\[0.50mm]
\hline
\end{tabular}
\end{table}
\begin{table}
\fbox{
 \epsfig{file=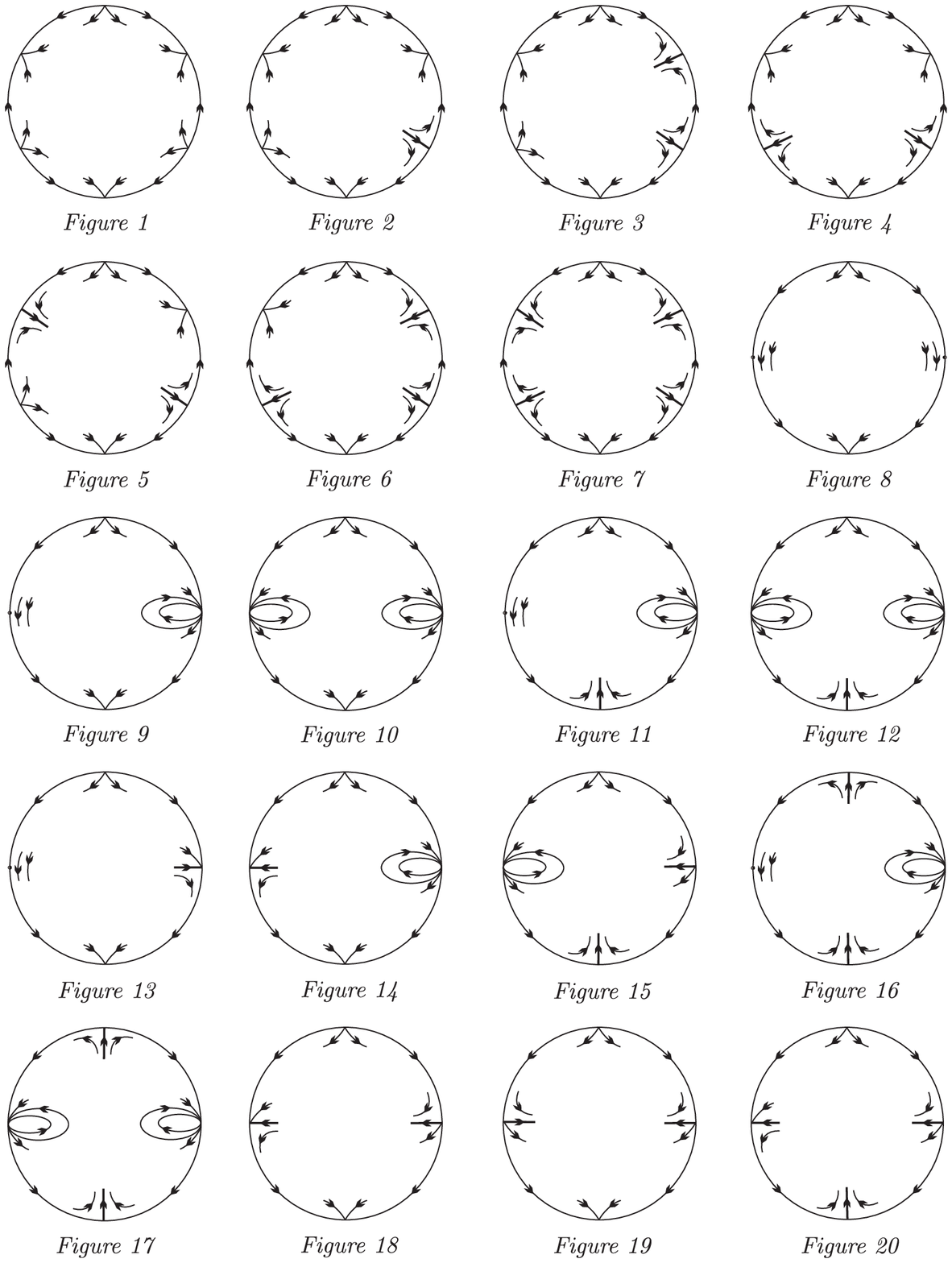,width=14.8cm}
}
\vspace{1mm}

\hfil Table 5 \hfill

\end{table}

\begin{table}
\fbox{
 \epsfig{file=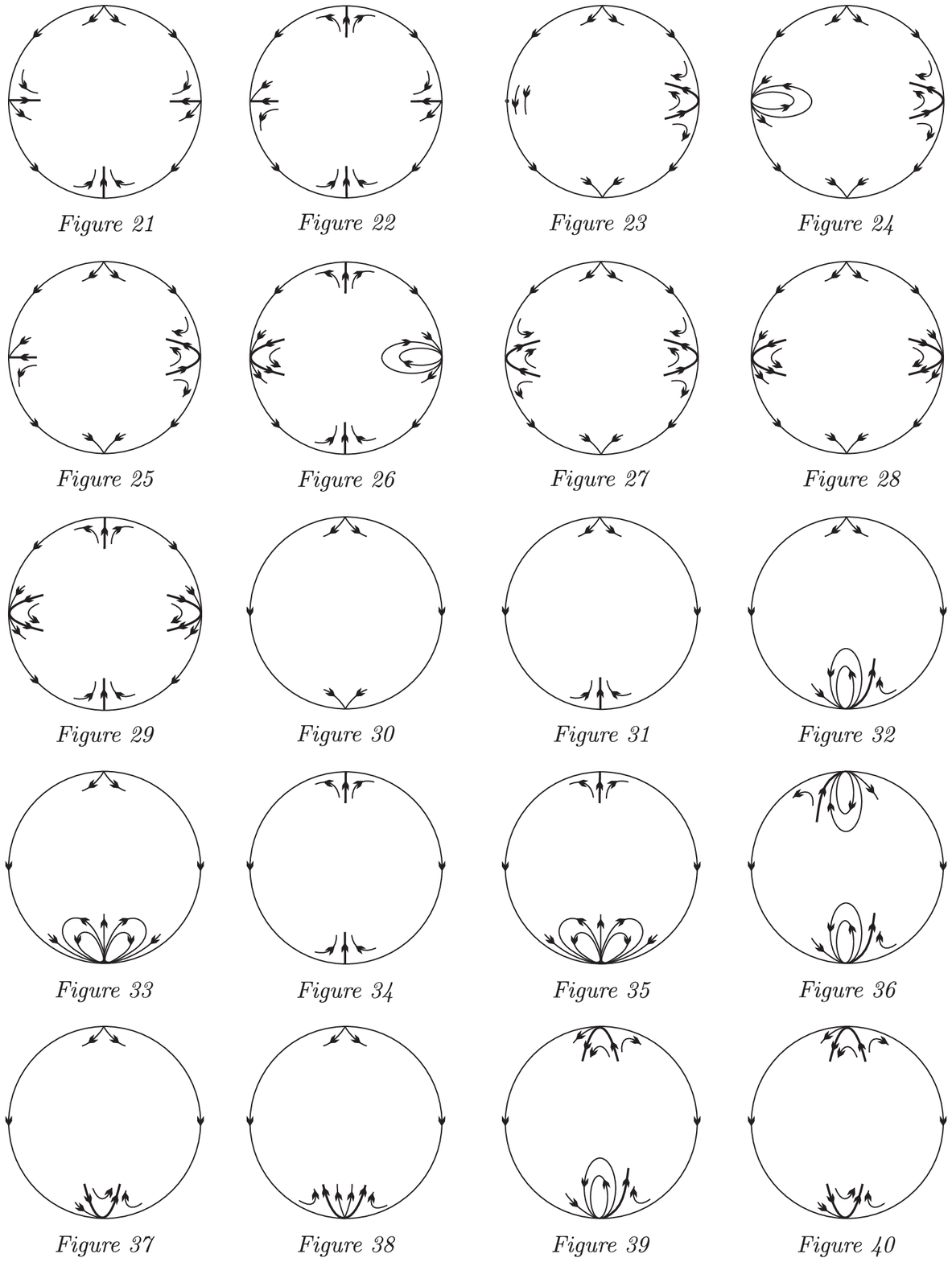,width=14.8cm}
}
\vspace{1mm}

\hfil Table 5 ({\it continued})  \hfill

\end{table}

 The proof is based on the Theorem
\ref{th_1} as well as on the invariant classification of quadratic systems at
infinity given in \cite{Nik_Vlp}, subject to some corrections as we shall
indicate below.

We point out that the affinely invariant conditions occurring  in part {\bf B}
of the theorem, greatly simplify the analogous conditions in \cite{Nik_Vlp}.

\brm[Corrections to {\rm \cite{Nik_Vlp}}]\label{rem:Errors}
In the statement of Theorem 2 ( a), b)) in \cite{Bautin}  (see page 92)
$\Delta_m>0$ must be replaced by $\Delta_m<0$ and  conversely. Since this
theorem was used in  \cite{Nik_Vlp} we have to note  that several expressions
in the sequences of the invariant conditions given in  \cite{Nik_Vlp} must be
taken with opposite sign, more precisely:
\vspace{-2mm}
\begin{itemize}
\vspace{-2mm}
    \item  {\sc Fig. 4}: the inequality $FS_1>0$ must be replaced by  $FS_1<0$;
\vspace{-2mm}
    \item  {\sc Fig. 5}: the inequality $FS_1<0$ must be replaced by $FS_1>0$;
\vspace{-2mm}
    \item  {\sc Fig. 6}: the inequality $GA<0$ must be replaced by $GA>0$;
\vspace{-2mm}
    \item  {\sc Fig. 7}: the inequality $GA>0$ must be  replaced by $GA<0$;
\vspace{-2mm}
    \item  {\sc Fig. 37}: the inequality $S_3<0$ must be replaced by $S_3<0,\ FS_1<0$;
\vspace{-2mm}
    \item  {\sc Fig. 38}: the inequalities $S_3>0,\ FS_1<0 $ must be replaced by $S_3>0,\ FS_1<0.$
\vspace{-3mm}
\end{itemize}
Furthermore  the saddle-node  given in {\sc Fig. 29} of \cite{Nik_Vlp}  is not
correctly placed.  The correct phase portrait is given here in Figure 15.
\erm

 {\it Proof of the Theorem} \ref{th_1}. {\bf A.} The phase portraits in
the vicinity of infinity of {\bf QS}$_{\bf ess}$ where obtained in
\cite{Nik_Vlp}.
 All calculations were done again for this article and as we
indicated in Remark \ref{rem:Errors}, all phase portraits obtained in
\cite{Nik_Vlp} with exception of {\sc Fig.} 29 turned out to be correct.
 Figure 29 in \cite{Nik_Vlp} needed to be modified at one of its singularities and we
give the respective corrected figure  in Table  5 ({\it Figure
15}).

In \cite{Nik_Vlp} (see pages 481--484), the phase portraits
appeared as they were obtained from calculations and not listed
according to their geometry. To draw attention to the geometry we
list them here for each possible value of $N_{\R}(S)$ according to
their topological complexity. In Table 3 we first place the number
$N_{\R}(S)$ of real singularities of the real foliation on
$\PP^2(\R)$, followed by the maximum number $\max(n_{\it sect})$
of sectors of singularities. Although these numbers could be read
on the value of $O(S)$, we place them in separate columns as they
are important invariants for the geometry at infinity of the
systems. We complete the table going through all phase portraits
and listing $O(S)$ which by itself determines uniquely 27 of the
40 phase portraits. To distinguish the remaining 13 phase
portraits we use the invariant $N_{\it hsect}^{f\infty
a}=(\,N_{\it hsect}^{f-a},\, N_{\it hsect}^{\infty-a}\,)$ whose
values we place in the last column, thus completing the
classification.

 {\bf B.}  As in the proof of part {\bf A} we  use the results
in \cite{Nik_Vlp} subject to  the  modifications in Remark \ref{rem:Errors}.
Since some letters appear both here and in \cite{Nik_Vlp}  but  not always with
the same meaning, we shall use the convention to apply   "tilde" to letters
which are used to denote comitants in \cite{Nik_Vlp}.

The proof of part {\bf B} proceeds in 3 steps:

I) In this step we replace the conditions in \cite{Nik_Vlp} subject to  the
modifications in Remark \ref{rem:Errors} with conditions involving newly
defined comitants and invariants as we shall indicate below.

II) In this part we simplify the conditions obtained in step I) in
order to obtain the corresponding conditions in the last column of
Table 4.

III) We prove that these last conditions are affinely invariant.

{\it Proof of step I}.   First of all we shall prove that the comitants used
in \cite{Nik_Vlp} (see Appendix) can be replaced respectively by the comitants
used here as follows:
\be\label{Corresp:1}
\begin{gathered}
\tilde\mu \Rightarrow \mu_0; \quad \tilde H \Rightarrow \mu_1; \quad \tilde G
\Rightarrow \mu_2; \quad  \tilde F \Rightarrow \mu_3; \quad \tilde V\Rightarrow
\mu_4; \quad \tilde L\Rightarrow C_2;\quad  \quad  \tilde M \Rightarrow M;\\
\quad  \tilde\eta \Rightarrow \eta;\quad
 \tilde \theta\Rightarrow \kappa; \quad
\tilde N\Rightarrow  K; \quad
 \tilde S_1\Rightarrow K_1;
 \quad  \tilde A \Rightarrow L;
 \quad \tilde A+4\tilde N\Rightarrow R;  \\
\tilde A+\tilde N\Rightarrow N; \quad  \tilde\sigma\Rightarrow \kappa_1; \quad
\tilde S_2\Rightarrow K_2; \quad \tilde S_3\Rightarrow K_3;
 \quad \tilde S_4\Rightarrow \kappa_2;
\end{gathered}
\ee
Indeed, firstly  the following relations among the comitants (\ref{Corresp:1})
hold:
\be\label{com_gen}
\begin{gathered}
   \mu_0=\tilde \mu, \quad
   \mu_1=2\tilde H, \quad
   \mu_2=\tilde G, \quad
   \mu_3=\tilde F, \quad
   \mu_4=\tilde V,\quad
    C_2=\tilde L,\quad
    M=8\tilde M,\\
    \eta=\tilde \eta,\quad
  \kappa= 64\tilde \theta, \quad
    K=4\tilde N, \quad
    K_1=\tilde S_1, \quad
    L=8\tilde A,\quad R=8(\tilde A + 4\tilde N).
\end{gathered}
\ee
Therefore we  only have to compare the conditions involving the comitants
\be\label{Compar:1}
N,\quad  \kappa_1,\quad \kappa_2,\quad K_2,\quad  K_3
\ee
and show the corresponding equivalence with the conditions involving  the
comitants
\be\label{Compar:2}
\tilde A+\tilde N,\quad  \tilde \sigma,\quad \tilde S_4,\quad \tilde S_2,\quad
\tilde S_3
\ee
  in \cite{Nik_Vlp}, respectively.

We point out that all comitants (\ref{Compar:2}) are only used  for the systems
$(\SSS_{I\!I\!I})$ and $(\SSS_{I\!V})$.  So, in what follows we shall  examine
each one of this cases.

We first consider the systems  of the form  $(\SSS_{I\!I\!I})$.

In this case we have four singularities on the equator  (i.e. $\eta=0,\
M\ne0$). The phase portraits in the  vicinity of infinity of these systems are
given by one of the Figures 8-29 both here and in \cite{Nik_Vlp}.
 One can observe, that all
comitants (\ref{Compar:1}) (respectively, (\ref{Compar:2})) are used for
systems $(\SSS_{I\!I\!I})$ only in the case when $\kappa=0$ (respectively,
$\tilde \theta=0$). In this case for systems $(\SSS_{I\!I\!I})$ the condition
$\kappa=-64h^2=0$ yields $h=0$ and we obtain the systems
\be \label{4s1}
  \dot x=k +cx +dy+gx^2,\qquad  \dot y=l+ex+fy +(g-1)xy,
\ee
for which $L=8gx^2$ and
$$
   \kappa_1=-32\,d,\ \  \tilde \sigma=-\frac{d}{4}(5g^2-2g+1);\quad
      N=(g-1)(g+1)x^2,\ \ \tilde A+\tilde N=\frac{1}{2}g(g+1)x^2.
$$
Clearly,  the condition $\kappa_1=0$ is equivalent to
$\tilde\sigma=0$. We now compare the signs of $N$ and $\tilde
A+\tilde N$.  As in Table 4 the comitant $N$ appears only in two
cases (i.e. Figures  16 and 26) and in these cases the condition
$L<0$ (i.e. $g<0$) is used, from the expressions of $N$ and $
\tilde A+\tilde N$ above we obtain   $\sign(N) =\sign(\tilde
A+\tilde N)$.

We observe  from Table 4  that the comitant $K_2$ is applied for
systems $(\SSS_{I\!I\!I})$ only when  $\kappa=\kappa_1=0$,
$L\ne0$. Since $\kappa_1=0$ implies  $d=0$ the systems (\ref{4s1})
become
\be \label{4s1a}
  \dot x=k +cx +dy+gx^2,\qquad  \dot y=l+ex+fy +(g-1)xy,
\ee
and we calculate\quad $
   K_2=48(g^2-g+2)(c^2-4gk)x^2, \ \ \tilde S_2=2g^2(c^2-4gk)x^2.
$\ Hence, $K_2$ has a well determined sign and since  for every $ g$ we have
$g^2-g+2>0$, from   $L\ne0$ we obtain $\sign(K_2)=\sign(\tilde S_2)$.

  We note that
the invariant $\kappa_2(a)$ is here used only to distinguish
Figures 8 and 17 in the case when systems ($\SSS_{I\!I\!I}$)
belong to the class $\Sigma_{33}$ in Table 2. Since for this class
the conditions $\kappa=L=K_1=0$ hold for systems
($\SSS_{I\!I\!I}$), we obtain respectively $h=g=c^2+d^2=0$. So,
the systems ($\SSS_{I\!I\!I}$) become
\be\label{4s1b}
  \dot x=k,\qquad  \dot y=l+ex+fy -xy,
\ee
for which $\kappa_2=-k, \quad \tilde S_4=-2k $ and, hence,
$\sign(\kappa_2)=\sign(\tilde S_4)$.

It remains to consider systems of the form $(\SSS_{I\!V})$.  For the Figures
30-40 which can occur for this class of systems, only the comitants $\tilde
S_2$ and $\tilde S_3$ of (\ref{Compar:2})  were used in \cite{Nik_Vlp}. Hence
we only have to examine the conditions given in terms of  comitants $K_2$ and
$K_3$ from (\ref{Compar:1}).

We observe that the comitant $K_2$ is used to distinguish Figures
30 and 32 when we also have $K\ne0$. In this case the systems
($\SSS_{I\!V}$) belong to the class $\Sigma_{17}$ in Table 2 with
conditions $\mu_0=\mu_1=0$.  For systems $(\SSS_{I\!V})$ we have $
\mu_0=-8h^3$. Hence $h=0$ and the systems ($\SSS_{I\!V}$) become
\be\label{5s4}
 \dot x= k+cx+dy+2gx^2,\quad
 \dot y=l+ex+fy-x^2+2gxy,
\ee
for which\ $ K=2g^2x^2,\  \mu_1=8dg^3x. $\ As $K\ne0$, the condition $\mu_1=0$
implies $d=0$ and  we obtain  the systems
\be\label{5s4a}
 \dot x= k+cx+dy+2gx^2,\quad
 \dot y=l+ex+fy-x^2+2gxy,
 \ee
for which we have:\quad
$
  K_2=24g^2(c^2-8gk)x^2, \ \ \tilde S_2=4g^2(c^2-8gk)x^2. \ \
$ Thus, in the case under consideration the comitant $K_2$ has a well
determined sign  and  $\sign(K_2)=\sign(\tilde S_2)$.

We examine now  the comitant $K_3$ which is applied for systems ($\SSS_{I\!V}$)
only in the cases when  $\Delta_S\ge3$, i.e. $\mu_{0,1,2}=0$. So, we shall
consider the systems (\ref{5s4})
 for which $\mu_0=0$ and we examine two subcases: $K\ne0$ and $K=0$.

If $K\ne0$ then  $g\ne0$ and for the systems~(\ref{5s4}) the condition
$\mu_1=0$ gives $d=0$. Moreover we may assume $e=f=0$ via a translation.
 So, we obtain the systems
\be\label{2s}
  \dot x=k +cx + 2gx^2,\qquad  \dot y=l - x^2 +2gxy,
\ee
for which $\mu_2=8g^3kx^2$ and as $g\ne0$ the condition $\mu_2=0$ yields $k=0$.
Then for the systems (\ref{2s}) we obtain $  K_3=-12g^2lx^6$, \ $ \tilde
S_3=-12g^2lx^6$.  Hence $K_3$ has a well determined sign and
$\sign(K_3)=\sign(\tilde S_3)$.

Assume now  $K=0$, i.e.   $g=0$ and for the systems (\ref{5s4}) we obtain
$\mu_1=0$, $\mu_2=d^2x^2$. Thus, the condition $\mu_2=0$ yields $d=0$ and we
obtain the following systems
\be\label{sys:K3}
 \dot x= k+cx,\qquad
 \dot y=l+ex+fy-x^2,
\ee
for which  $K_3=3f(2c-f)x^6=\tilde S_3$. \EProof

{\it Proof of step II}. We show below how some  of the conditions
in \cite{Nik_Vlp}  can be substituted by simpler ones  in Table 4.
To do this we shall prove the following five lemmas.

\blm\label{equiv:17-22} Let\ $
\tilde {\frak C}$ be the conjunction of the all the conditions: $\tilde
\eta=\tilde \mu=\tilde H=\tilde \theta=\tilde\sigma=0$\ and $\tilde M\tilde
G\tilde A\ne0 $. Let $ {\frak C}$ be the conjunction of the following
conditions: $ \eta=\mu_0=\mu_1=\kappa=\kappa_1=0\ \ \mbox{and}\ \ M\mu_2L\ne0.
$ We have the following equivalences:
$$
\ba{llll}
Figure\ 8: &\tilde {\frak C},\,  \tilde G\ne0,\, \tilde A>0,\, \tilde S_2<0 &
\Leftrightarrow  &{\frak C},\,
\mu_2>0,\, L>0,\ K_2\!<\!0;\\
Figure\ 10: &  \tilde {\frak C},\,\tilde G<0,\, \tilde A>0,\, \tilde S_2>0,\
\tilde N<0 & \Leftrightarrow  &{\frak C},\,
\mu_2<0,\, L>0,\, K<0;\\
Figure\ 17: & \tilde {\frak C},\, \tilde G\ne0,\, \tilde A<0,\, (\tilde
S_2\!\le\!0)\lor (\tilde G\!>\!0,\, \tilde S_2\!>\!0) &  \Leftrightarrow
&{\frak C},\,
\mu_2>0,\, L<0;\\
Figure\ 19: &\tilde {\frak C},\,  \tilde G\ne0,\, \tilde A>0,\, (\tilde
S_2\!=\!0)\lor (\tilde G\!>\!0,\, \tilde S_2\!>\!0) &  \Leftrightarrow  &{\frak
C},\,
\mu_2>0,\, L>0,\, K_2\!\ge\!0;\\
Figure\ 27: &\tilde {\frak C},\,  \tilde G<0,\, \tilde A>0,\, \tilde S_2>0,\,
\tilde N>0 & \Leftrightarrow  &{\frak C},\,
\mu_2<0,\, L>0,\, K>0;\\
Figure\ 29: & \tilde {\frak C},\, \tilde G<0,\, \tilde A<0,\, \tilde S_2>0 &
\Leftrightarrow &
 {\frak C},\,\mu_2<0,\, L<0.\\
\ea
$$
\elm
\BProof According to \eqref{com_gen} the conditions $\tilde {\frak C}$ and $
{\frak C}$ are equivalent. We are in the class of the systems
$(\SSS_{I\!I\!I})$ for which we must apply the conditions on the
right, i.e. $\mu_0=\mu_1=0$, $\mu_2\ne0$, and $\kappa=\kappa_1=0$,
$L\ne0$. For the systems $(\SSS_{I\!I\!I})$ we have $\
\kappa=-64h^2$, $\kappa_1=-32d$ and hence conditions
$\kappa=\kappa_1=0$ yield $h=d=0$. Then
$$
\mu_0=\mu_1=0,\quad \mu_2=g[f^2g+cf(g-1) +k(g-1)^2]x^2\ne0
$$
and since $g\ne0$ we may assume $c=0$ via a translation. Hence we get the
systems
\be\label{Sys:S32a}
  \dot x=k+gx^2,\quad
  \dot y=l+ex+fy+ (g-1)xy,
\ee
for which
\be\label{values3:1a}
\bal
&\mu_{0,1}=0,\ \mu_2=g[f^2g+k(g-1)^2]x^2\tilde G\ne0,\ L=gx^2=8\tilde A\ne0,\\
&K=2g(g-1)x^2=4\tilde N,\ K_2=-192gk(g^2-g+2)x^2,\ \tilde S_2=-8g^3k.
\eal
\ee
We observe, that  $\sign(K_2)=\sign(\tilde S_2)$ because the discriminant of
the quadratic polynomial $g^2-g+2$ is negative. We shall consider two cases:
$L<0$ and $L>0$.

\noindent{\bf Case $L<0$.} If $\mu_2<0$ (then $\tilde G<0$)  from
(\ref{values3:1a}) it follows that $\tilde S_2>0$ and hence we obtain the
conditions indicates on the left in the lemma, which correspond to Figure 29.
Thus the conditions $L<0$ and $\mu_2<0$ lead to Figure 29.

Assume $\mu_2>0$ (then $\tilde G>0$).  If either  $K_2>0$ (then $\tilde S_2>0$)
or $K_2\le0$ (then $\tilde S_2\le0$) we obtain the conditions on the left for
Figure 17. Taking into account that for $\mu_2\ne0$ from (\ref{values3:1a}) it
follows that the condition $\tilde S_2\le0$ implies $\mu_2>0$ (then $\tilde
G>0$) we conclude, that the conditions $L<0$ and $\mu_2>0$ lead to Figure 17.

\noindent{\bf Case $L>0$.} Suppose firstly  $\mu_2<0$.  Then $\tilde G<0$  and
from  (\ref{values3:1a}) we have $\tilde S_2>0$ and  $\tilde N\ne0$ (i.e.
$K\ne0$). Hence we obtain the conditions for Figure 10 (on the left in the
lemma)  if $K<0$ and for Figure 27 if $K>0$.

Assume now $\mu_2>0$ (then $\tilde G>0$).  From (\ref{values3:1a}) we obtain
$\tilde S_2\le0$ (then $K_2\le0$) which yields $\mu_2>0$. Hence we conclude,
that the conditions $L>0,$ $\mu_2>0$, $K_2\ge0$ lead to the Figure 19, whereas
the conditions $L>0,$ $\mu_2>0$, $K_2<0$ lead to the Figure 8. \EProof

\blm\label{equiv:8,20,}
Let\ $ \tilde {\frak C}_1 $ be the conjunction of the following conditions: $
\tilde \eta=\tilde \mu=\tilde H=\tilde G=\tilde F=\tilde \theta=\tilde\sigma=0
\ \mbox{and} \ \tilde M\tilde V\tilde A\ne0 $. Let $ {\frak C}_1 $ be the
conjunction of the following conditions: $
\eta=\mu_0=\mu_1=\mu_2=\mu_3=\kappa=\kappa_1=0\ \mbox{and} \ M\mu_4L\ne0. $ We
have the following equivalences:
$$
\ba{llll}
Figure\ 8: &\tilde {\frak C}_1,\,  \tilde V\ne0,\, \tilde A\ne0,\, \tilde
N=0,\, \tilde S_2<0 &\ \Leftrightarrow\ &{\frak C}_1,\,
\mu_4>0,\, L>0,\, K_2<0;\\
Figure\ 17: &\tilde {\frak C}_1,\,  \tilde V\ne0,\, \tilde N\ne0,\, \tilde A<0
&\ \Leftrightarrow\ &{\frak C}_1,\,
\mu_4>0,\, L<0;\\
Figure\  18:\  &  \left[\!\!\ba{l}\tilde {\frak C}_1,\, \tilde A \tilde V\ne0,\, (\tilde N=0,\, \tilde S_2=0)\\
                            \lor(\tilde N\ne0,\, \tilde A>0,\, \tilde A+4\tilde
                            N\ge0)\\
                            \lor(\tilde N=0,\, \tilde V>0,\, \tilde S_2>0)\ea\!\!\right]
                            & \Leftrightarrow &\!\!
\left[\!\!\ba{l}{\frak C}_1,\, \mu_4>0,\, L>0,\\ ( R\ge0,\, K\ne0)\lor\\
      (K_2\ge0,\, K=0)\ea\!\!\right];\\[7mm]
Figure\ 24: &\tilde {\frak C}_1,\,   \tilde N=0,\, \tilde A\ne0,\, \tilde V<0
&\ \Leftrightarrow\ &{\frak C}_1,\,
\mu_4<0,\, L\ne0;\\
Figure\ 28: &\tilde {\frak C}_1,\,  \tilde V \tilde N\ne0,\, \tilde A>0,\,
\tilde A+4\tilde N<0 &\ \Leftrightarrow\ &{\frak C}_1,\,
\mu_4>0,\, L>0,\, R<0.\\
\ea
$$
\elm
\BProof We are in the class of systems $(\SSS_{I\!I\!I})$ for which we must set
the conditions $\mu_0=\mu_1=\mu_2=\mu_3=0$, $\mu_4\ne0$, and
$\kappa=\kappa_1=0$, $L=8\tilde A\ne0$. It was shown before (see page
\pageref{Sys:S32a}) that for the systems $(\SSS_{I\!I\!I})$ the conditions
$\kappa=\kappa_1=0$ yield $h=d=0$. Then $L=gx^2\ne0$ and $K=2g(g-1)x^2$ and we
shall construct two canonical forms corresponding to the cases $K\ne0$ and
$K=0$.

 Assume firstly  $K\ne0$. Then $g-1\ne0$ and we may assume $e=f=0$ due to a
translation.
 Therefore considering the  conditions $h=d=e=f=0$, for the systems  $(\SSS_{I\!I\!I})$
calculations yield:  $\mu_0=\mu_1=0$, $\mu_2=gk(g-1)^2$ and by $g(g-1)\ne0$ the
condition $\mu_2=0$  yields $k=0$. This implies $\mu_3=-clg(g-1)x^3$,
$\mu_4=lx^3[lg^2x+c^2(g-1)y]$. Hence, the conditions $\mu_3=0$ and $\mu_4\ne0$
yield $c=0$ and we get the systems
\be\label{Sys:S31}
  \dot x=gx^2,\quad
  \dot y=l+(g-1)xy,
\ee
for which
\be\label{values3:1}
\bal
&\mu_{0,1,2,3}=0,\ \mu_4=g^2l^2x^4=\tilde V,\ L=8gx^2=8\tilde A\ne0,\
K_2=0=\tilde S_2,\\
&K=2g(g-1)x^2=4\tilde N\ne0,\ R=8g(2g-1)x^2=8(\tilde A+4 \tilde N).
\eal
\ee

Suppose now that the condition  $K=2g(g-1)x^2=0$ holds.  Since $L=gx^2\ne0$
this yields $g=1$ and we may assume $c=0$ via a translation. Then we obtain
$\mu_2=f^2x^2=0$ which implies $f=0$ and we get the systems
\be\label{Sys:S32}
  \dot x=k+x^2,\quad
 \dot y=l+ex,
\ee
for which
\be\label{values3:2}
\bal
& \mu_{0,1,2,3}=0,\ \mu_4=(l^2+ke^2)x^4=\tilde V\ne0,\ L=8x^2=8\tilde A,\\
& K=0=\tilde N, \ R=8x^2=8(\tilde A+4 \tilde N),\ K_2=-384\,kx^2=48\tilde S_2.
\eal
\ee
We shall consider two cases: $\mu_4<0$ and $\mu_4>0$.

\noindent{\bf Case $\mu_4<0$.} Then $\tilde V<0$ and  from
(\ref{values3:1}) and (\ref{values3:2}) we have the conditions  $\tilde N=0$
and $\tilde S_2>0.$ Hence the conditions  $\mu_4<0$ and $L\ne0$ lead to the
conditions in the lemma corresponding to Figure 24.

\noindent{\bf Case $\mu_4>0$.} In this case $\tilde V>0$ and we shall examine
two subcases: $L<0$ and $L>0$.

{\bf Subcase $L<0$.}  Then $\tilde A<0$. From (\ref{values3:1}) and
(\ref{values3:2}) we conclude that  $\tilde N\ne0$ and we obtain the conditions
corresponding to  Figure 17. Hence we conclude that for $\mu_4>0$ and $L<0$ we
get Figure 17.

{\bf Subcase $L>0$.}  Hence $\tilde A>0$.

{\it a)} If  $R<0$ (then $\tilde A+4 \tilde N<0$) from (\ref{values3:1}) and
(\ref{values3:2}) we obtain  $\tilde N\ne0$ and hence we get the conditions for
Figure 28.

{\it b)} Assume now $R\ge0$. If $K\ne0$ (then $\tilde N\ne0$)  we obtain one
sequence of conditions for Figure 18, and namely: $\tilde N\ne0$, $\tilde A>0$
and $\tilde A+4 \tilde N\ge0$.

Suppose $K=0$ (i.e. $\tilde N=0$). If in addition  $K_2<0$ (then $\tilde
S_2<0$) then we obtain the conditions for Figure 8. From (\ref{values3:1}) and
(\ref{values3:2}) we obtain that the condition $K_2<0$ implies $\tilde N=0$.
Then we conclude, that for $\mu>0,$ $L>0$ and $K_2<0$ we obtain the conditions
for Figure 8.

Assuming $K_2\ge0$ (then $\tilde S_2\ge0$) and taking into account that we are
in the case $\mu_4>0$, we get two of the series of conditions for Figure 18,
which can be combined into the following series: $\mu_4>0,$ $ K=0,$ $ L>0,$ $
K_2\ge0$ . \EProof

\blm\label{equiv:30,34,35}
Let\ $ \tilde {\frak C}_2 $ be the conjunction of all the conditions: $ \tilde
M=\tilde \mu=\tilde H=0 $ and $ \tilde L\tilde G\ne0 $. Let $ {\frak C}_2 $ be
the conjunction of the following conditions: $ M=\mu_0=\mu_1=0\ \mbox{and} \
C_2\mu_2\ne0. $ We have the following equivalences:
$$
\!\!\ba{llll} Figure\ 30:  &\! \tilde {\frak C}_2, \tilde G\!\ne\!0, (\tilde
N\!\ne\!0, \tilde S_2<0)\lor (\tilde N\!=\!0) & \Leftrightarrow &{\frak C}_2,
\mu_2\!>\!0, (K\!\ne\!0, K_2\!<\!0)\lor (K\!=\!0);\\
Figure\ 32:  &\!\tilde {\frak C}_2,  \tilde G\tilde N\!\ne\!0, (\tilde G\!>\!0,
\tilde S_2\!>\!0)\lor (\tilde S_2\!=\!0)\!\! &\Leftrightarrow &{\frak C}_2,
\mu_2>0, K\ne0, K_2\ge0;\\
Figure\ 40:  &\!\tilde {\frak C}_2,  \tilde G<0,\ \tilde N\ne0,\ \tilde S_2>0 &
\Leftrightarrow &{\frak C}_2,
\mu_2<0.\\
\ea
$$
\elm
\BProof We are in the class of systems $(\SSS_{I\!V})$ for which we must set
the conditions $\mu_0=\mu_1=0$, $\mu_2\ne0$. We have $\mu_0=-h^3=0$ which
implies $ h=0$ and then $\mu_1=dg^3x$ and $K=2gx^2$. We shall consider two
subcases: $K\ne0$ and $K=0$.

Assume firstly $K\ne0$. Then $g\ne0$  and the condition $\mu_1=0$ yields $d=0$.
We can assume $g=1$ and $e=f=0$ due to the rescaling $x\to x/g$, $y\to y/g^2$
and
 a translation. Then  we get the systems
\be\label{Sys:S41a}
  \dot x =k+cx+x^2,\quad
  \dot y=l-x^2+xy,
\ee
for which
\be\label{values:1a}
\mu_{0,1}=0,\ \mu_2=kx^2=\tilde G\ne0,\ K=2x^2=4\tilde N,\
K_2=48(c^2-4k)x^2=24\tilde S_2.
\ee

Admit now $K=0$. Hence $g=0$ and  we can  assume $e=0$ due to a translation.
Then we   obtain the systems
\be\label{Sys:S42a}
  \dot x=k+cx+dy,\quad
  \dot y=l+fy-x^2,
\ee
for which
\be\label{values:2a}
\mu_{0,1}=0,\ \mu_2=d^2x^2=\tilde G\ne0,\ K=0=\tilde N,\ L_2=0=\tilde S_2.
\ee
\noindent{\bf Case $\mu_2<0$}. From (\ref{values:1a}) and (\ref{values:2a})
it follows that the condition $\mu_2<0$ implies   $\tilde N\ne0$ and $\tilde
S_2>0$. Hence we obtain the  conditions for Figure 40 and we conclude that the
condition $\mu_2<0$ immediately leads to the conditions for  Figure 40.

\noindent{\bf Case $\mu_2>0$} ( i.e. $\tilde G>0$). Assume that the condition
$K\ne0$ holds (then $\tilde N\ne0$). If $K_2<0$ we have $\tilde S_2<0$ and then
we obtain the conditions for Figure 30. If either $K_2>0$ or  $K_2=0$ via
$\tilde G>0$ in both cases we  get Figure 32.

Suppose   $K=0$ (i.e. $\tilde N=0$).  In this case we  obtain the conditions
$\tilde G\ne0$, $\tilde N=0$ which lead to  Figure 30. Note that from
(\ref{values:1a}) and (\ref{values:2a}) it follows that  the condition $K=0$
implies  $\mu_2>0$.
 \EProof

\blm\label{equiv:32,37,38}
Let\ $ \tilde {\frak C}_3 $ be the conjunction of the following conditions: $
\tilde M=\tilde \mu=\tilde H=\tilde G=\tilde N=0 \ \mbox{and} \ \tilde L\tilde
F\ne0$. Let $ {\frak C}_3$ be the conjunction of the following conditions: $
M=\mu_0=\mu_1=\mu_2=K=0\ \mbox{and} \ C_2\mu_3\ne0. $ We have the following
equivalences:
$$
\ba{llll}
Figure\ 31:\  & \left[\!\!\ba{l}\tilde {\frak C}_3,\, \tilde F\ne0,  (\tilde S_3=0)\\
 \lor(\tilde F\tilde S_1>0,\, \tilde S_3>0)\ea\!\!\right]  &\ \Leftrightarrow\ &{\frak C}_3,\,
 \mu_3K_1>0,\, K_3\ge0;\\[5mm]
Figure\ 33:\  &\tilde {\frak C}_3,\,   \tilde F\tilde S_1<0, \tilde S_3<0 &\
\Leftrightarrow\
& {\frak C}_3,\, \mu_3K_1<0;\\
Figure\ 38:\  &\tilde {\frak C}_3,\,   \tilde F\tilde S_1>0, \tilde S_3<0 &\
\Leftrightarrow\
& {\frak C}_3,\, \mu_3K_1>0,\, K_3<0.\\
\ea
$$
\elm
\BProof We are in the class of systems $(\SSS_{I\!V})$ for which we must set
the conditions $\mu_0=\mu_1=\mu_2=0=K$, $\mu_3\ne0$. We have
$\mu_0=-h^3=0$ hence $ h=0$ and then  $K=2gx^2$. The condition
$K=0$ yields $g=0$ and this leads to the systems (\ref{Sys:S42a})
for which the condition
 $\mu_2=d^2x^2=0$ yields $d=0$. Hence  we obtain the
systems
\be\label{Sys:S420}
  \dot x=k+cx,\quad
  \dot y=l+fy-x^2,
\ee
for which
\be\label{values:2ab}
\bal
&\mu_{0,1,2}=0,\quad \mu_3=-c^2fd^2x^3=\tilde F\ne0,\quad K=0=\tilde N,\\
&K_1=-cx^3=\tilde S_1,\quad\ K_3=6f(2c-f)x^6=\tilde S_3.
\eal
\ee
We note that $\mu_3K_1=c^3fx^6\ne0$ and hence
$\sign(\mu_3K_1)=\sign(cf)=\sign(\tilde F\tilde S_1)$.

\noindent{\bf Case  $\mu_3K_1<0$}.
From (\ref{values:2ab})  we obtain $\tilde S_3<0$ and hence we conclude that
the condition $\mu_3K_1<0$ leads to the conditions for Figure  33.

\noindent{\bf Case  $\mu_3K_1>0$.} For $K_3<0$ (then $\tilde
S_3<0$)  we obtain  the conditions for Figure 38. If either
$K_3>0$ or  $K_3=0$  we observe that in both cases we get the
conditions for Figure 31. From (\ref{values:2ab}) it follows that
the condition $K_3=0$ implies $\mu_3K_1>0$. Therefore  we conclude
that the conditions   $\mu_3K_1>0$ and $K_3\ge0$ lead to the
conditions for  Figure 31.
 \EProof

\blm\label{equiv:30,34,39,40}
Let\ $ \tilde {\frak C}_4 $ be the conjunction of the following conditions: $
\tilde M=\tilde \mu=\tilde H=\tilde G=\tilde F=0 \ \mbox{and} \ \tilde L\tilde
V\ne0$. Let $ {\frak C}_4$ be the conjunction of the following conditions: $
M=\mu_0=\mu_1=\mu_2=\mu_3=0\ \mbox{and} \ C_2\mu_4\ne0. $ We have the following
equivalences:
$$
\ba{llll}
Figure\ 30:\  & \left[\!\!\ba{l}\tilde {\frak C}_4,\, \tilde V\ne0,\, (\tilde N\ne0,\, \tilde S_3>0)\\
 \lor(\tilde N=\tilde S_1= \tilde S_3=0)\lor\\
  (\tilde N=0,\, \tilde S_1\ne0,\, \tilde V>0) \ea\!\!\right]  &\ \Leftrightarrow\ &{\frak C}_4,\,
 \mu_4>0,\, K_3\ge0;\\[7mm]
Figure\ 32:\  &\tilde {\frak C}_4,\, \tilde V\ne0,\, \tilde N=\tilde S_1=0,\,
\tilde S_3\ne0&\ \Leftrightarrow\
& {\frak C}_4,\, \mu_4>0,\, K_3<0,\, K=0;\\
Figure\ 35:\  & \tilde {\frak C}_4,\,\tilde V<0,\, \tilde N=0,\, \tilde
S_1\ne0&\ \Leftrightarrow\
& {\frak C}_4,\, \mu_4<0;\\
Figure\ 36:\  &\tilde {\frak C}_4,\, \tilde V\ne0,\, \tilde N\ne0,\, \tilde
S_3<0&\ \Leftrightarrow\
&{\frak C}_4,\,  \mu_4>0,\, K_3<0,\, K\ne0.\\
\ea
$$
\elm
\BProof We are in the class of systems $(\SSS_{I\!V})$ for which we must set
the conditions $\mu_0=\mu_1=\mu_2=\mu_3=0$, $\mu_4\ne0$. We have $\mu_0=-h^3=0$
which implies $ h=0$ and then $\mu_1=dg^3x$ and $K=2gx^2$. We shall consider
two subcases: $K\ne0$ and $K=0$.

If  $K\ne0$ then the  condition $\mu_1=0$  leads to the systems
(\ref{Sys:S41a}) for which  $\mu_2=kx^2.$ Hence the condition $\mu_2=0$ yields
 $k=0$ and we calculate: $\mu_3=-clx^3$ and $\mu_4=-l(c^2x-lx-c^2y)x^3$. Hence
the conditions $\mu_3=0$ and $\mu_4\ne0$ yield $c=0$, $l\ne0$ and we obtain the
systems
\be\label{Sys:S41oa}
  \dot x =x^2,\quad
  \dot y=l-x^2+xy,
\ee
for which
\be\label{values:1}
\mu_{0,1,2,3}=0,\ \mu_4=l^2x^4=\tilde V\ne0,\ K=\frac{1}{2}x^2=4\tilde N,\
 K_3=-6lx^6=\tilde S_3\ne0.
\ee

 Admit now that $K=0$. This leads to the systems (\ref{Sys:S420}) for which
the condition $\mu_3=-c^2fd^2x^3=0$  yields $cf=0$. Then  we get the systems
\be\label{Sys:S42oa}
  \dot x=k+cx,\quad
  \dot y=l+fy-x^2,
\ee
with $cf=0$ and
\be\label{values:2}
\bal
&\mu_{0,1,2,3}=0,\quad \mu_4=(k^2-c^2l)x^4=\tilde V\ne0,\quad  K=0=\tilde N,\\
&K_1=-cx^3=\tilde S_1,\quad K_3=-6f^2x^6=\tilde S_3,\quad K_1K_3=0.
\eal
\ee
\noindent{\bf Case $\mu_4<0$} (i.e.  $\tilde V<0$).
From (\ref{values:1}) and  (\ref{values:2}) we obtain that the condition
$\mu_4<0$ implies  $\tilde N=0$ and $ \tilde S_1\ne0$. Hence for $\mu_4<0$ we
obtain the conditions for Figure 35.

\noindent{\bf Case $\mu_4>0$}. Then $\tilde N>0$ and we shall consider 3 subcases:
$K_3<0$, $K_3>0$ and $K_3=0$.\\
\indent{\bf Subcase $K_3<0$.}
If $K\ne0$ then $\tilde N\ne0$ and  we have the conditions for Figure 36.
Suppose $K=0$, i.e. $\tilde N=0$. Then by $K_3\ne0$ from   (\ref{values:2}) we
have $\tilde S_1=0$. Therefore we conclude that conditions $K_3<0$ and  $K=0$
lead to the Figure 32.\\
\indent{\bf Subcase $K_3>0$.} Then $\tilde S_3>0$ and from (\ref{values:1})
and  (\ref{values:2}) we conclude that $K\ne0$, i.e. $\tilde N\ne0$.  Hence we
obtain one series of the conditions for Figure
30. \\
\indent{\bf Subcase $K_3=0$.} Then $\tilde S_3=0$ and according to
(\ref{values:1}) and  (\ref{values:2}) we have $K=0$. This leads to the systems
(\ref{Sys:S42oa}) for which the condition $K_3=0$ yields  $f=0$. Then we have
either $K_1\ne0$ (i.e. $\tilde S_1\ne0$) or $K_1=0$ (i.e. $\tilde S_1=0$).
Since the conditions $\tilde V>0$ and $\tilde S_3=0$ hold,  both cases lead to
the conditions for  Figure 30.

Lemma \ref{equiv:30,34,39,40} is proved and this completes the proof of the
step II. \EProof

{\it Proof of step III}. We draw the attention to the fact that all the
constructed  polynomials which were  used in Theorems \ref{th_1} and \ref{th_2}
are $GL$-comitants. But in  fact we are interested in the action of the affine
group $\Aff(2,\mathbb R)$ on these systems. We shall prove the following lemma.

\blm\label{Table:Propreties}
   The  polynomials which are used
   in  Theorems \ref{th_1} or  \ref{th_2}  have the properties
indicated in the Table 6. In the last column are indicated the
algebraic sets on which the  $GL$-comitants on the left are
$CT$-comitants.
  The Table  6 shows us   that all conditions
included in the statements of  Theorems \ref{th_1} or \ref{th_2} are affinely
invariant.
\elm
\begin{table}[!htb]
\begin{center}
\begin{tabular}{|c|c|c|c|c|c|}
\multicolumn{6}{r}{\bf Table 6}\\[2mm]
\hline \raisebox{-0.7em}[0pt][0pt]{Case} &
\raisebox{-0.7em}[0pt][0pt]{$GL$-comitants}
 & \multicolumn{2}{c|}{Degree in }  & \raisebox{-0.7em}[0pt][0pt]{Weight} &  Algebraic subset \\
\cline{3-4}
  &  & $\ \ a\ \ $ & $\!x$ and $y\!$ &    &  $V(*)$   \\
\hline\hline \rule{0pt}{2ex}
 $1$ & $\eta(a)$,\ $\mu_0(a)$,\ $\kappa(a)$ & $4$ &  $0$ & $ 2$  &  $V(0)$ \\[0.5mm]
\hline\rule{0pt}{2ex} $2$  & $C_2(a,x,y)$   & $1$ &  $3$ & $-1$  & $V(0)$\\[0.5mm]
\hline\rule{0pt}{2ex} $3$ & $K(a,x,y)$  & $2$ &  $2$ & $ 0$   &  $V(0)$ \\[0.5mm]
\hline\rule{0pt}{2ex} $4$ & $L(a,x,y)$  & $2$ &  $2$ & $ 0$   &  $V(0)$ \\[0.5mm]
\hline\rule{0pt}{2ex} $5$ & $M(a,x,y)$   & $2$ &  $2$ & $ 0$   &  $V(0)$ \\[0.5mm]
\hline\rule{0pt}{2ex} $6$ & $N(a,x,y)$ & $2$ &  $2$ & $ 0$   &  $V(0)$ \\[0.5mm]
\hline\rule{0pt}{2ex} $7$ & $R(a,x,y)$  & $2$ &  $2$ & $ 0$   &  $V(0)$ \\[0.5mm]
\hline \rule{0pt}{2ex} $8$ & $\kappa_1(a)$  & $3$ &  $0$ & $1$   & $V(\eta,\kappa)$ \\[0.5mm]
\hline\rule{0pt}{2ex}  $9$ & $\kappa_2(a)$  & $2$ &  $0$ & $0$   &  $V(\eta,\kappa,L,K_1)$ \\[0.5mm]
\hline\rule{0pt}{2ex} $10$ & $K_2(a,x,y)$  & $4$ &  $2$ & $ 0$   &  $V(\eta,\mu_0,\mu_1,\kappa,\kappa_1)$ \\[0.5mm]
\hline\rule{0pt}{2ex} $11$ & $K_3(a,x,y)$  & $4$ &  $6$ & $-2$   &  $V(M,\mu_0,\mu_1,\mu_2)$ \\[0.5mm]
\hline\rule{0pt}{2ex} $12$ & $K_1(a,x,y)$ & $2$ &  $3$ & $-1$   &  $V(K)$ \\[0.5mm]
\hline\rule{0pt}{2ex} $13$ & $\mu_1(a,x,y)$ & $4$ &  $1$ & $ 1$   &  $V(\mu_0)$ \\[0.5mm]
\hline\rule{0pt}{2ex} $14$ & $\mu_2(a,x,y)$ & $4$ &  $2$ & $ 0$   &  $V(\mu_0,\mu_1)$ \\[0.5mm]
\hline\rule{0pt}{2ex} $15$ & $\mu_3(a,x,y)$  & $4$ &  $3$ & $-1$   &  $V(\mu_0,\mu_1,\mu_2)$ \\[0.5mm]
\hline\rule{0pt}{2ex} $16$ & $\mu_4(a,x,y)$  & $4$ &  $4$ & $-2$   &  $V(\mu_0,\mu_1,\mu_2,\mu_3)$ \\[0.5mm]
\hline
\end{tabular}
\end{center}
\end{table}
\BProof  {\it I. Cases  1--7}.  The polynomials $\eta(a),\ \kappa(a)$,
$\mu_0(a)$, $K(a,x,y)$, $L(a,x,y),$ $ M(a,x,y),$ $ N(a,x,y)$ and $R(a,x,y)$ are
$T$-comitants, because these $GL$-comitants were constructed only by using the
coefficients of the polynomials $p_2(x,y)$ and $q_2(x,y)$.

{\it II. Cases 8--11}. $a)$ We consider the $GL$-invariant
$\kappa_1(a)$ which according to Table 4 was used only  in the
class of systems $(\SSS_{I\!I\!I})$. It was shown before (see page
\pageref{4s1}) that for $\kappa=0$ the systems $(\SSS_{I\!I\!I})$
can be brought by an affine transformation to the systems
(\ref{4s1}) for which $\kappa_1=-32d$.  On the other hand  for any
system in the orbit under the translation group action of a system
(\ref{4s1}) corresponding to a point $\ab\in \R^{12}$   we obtain
$\kappa_1(\ab)=-32\,d$. Hence the value of $\kappa_1$ does not
depend of the vector defining the  translations. Therefore we
conclude that
 the polynomial $\kappa_1$  is a $CT$-comitant modulo $\langle\eta,\kappa\rangle$.

$b)$ We consider now the $GL$-invariant $\kappa_2(a)$. From Table
4 we observe that $\kappa_2(a)$ is only applied to distinguish the
Figures 8 and 17 when for the  systems $(\SSS_{I\!I\!I})$ the
conditions $\kappa=L=K_1=0$ hold. As it was shown before (see page
\pageref{4s1b})  for $\kappa=L=K_1=0$ the systems
$(\SSS_{I\!I\!I})$ can be brought by an affine transformation to
the systems (\ref{4s1b}) for which $\kappa_2=-k$.  On the other
hand  for any system  in the orbit under the translation group
action of a system (\ref{4s1b}) corresponding to a point $\ab\in
\R^{12}$
 we obtain $\kappa_2(\ab)=-k$. Hence  we conclude
that the polynomial $\kappa_2$  is a $CT$-comitant modulo
$\langle\eta,\kappa,L,K_1\rangle$.

$c)$ We examine now the  $GL$-invariant $K_2(a)$ which was used in  cases
 $(\SSS_{I\!I\!I})$ and $(\SSS_{I\!V})$.  Assume firstly  $\eta=0$ and
 $M\ne0$, i.e. we are in the class of the systems $(\SSS_{I\!I\!I})$.
We have shown before (see page \pageref{4s1a}) that  for $\kappa=\kappa_1=0$
the systems $(\SSS_{I\!I\!I})$ can be brought by an affine transformation to
the systems (\ref{4s1a}) for which $K_2=48(g^2-g+2)(c^2-4gk)x^2$. Suppose  now
that the conditions
 $M=0$ and $C_2\ne0$ hold, i.e. we are in the class of the systems $(\SSS_{I\!V})$.
It was shown before (see page \pageref{5s4a}) that  for $\mu_0=\mu_1=0$ the
systems $(\SSS_{I\!V})$ can be brought by an affine transformation to the
systems (\ref{5s4a}) for which $K_2=24g^2(c^2-8gk)x^2$.

On the other hand for any system  in the orbit under the translation group
action of a system (\ref{4s1a}) (respectively, of a system (\ref{5s4a}))
corresponding to a point $\ab\in \R^{12}$ (respectively, $\ab_1\in \R^{12}$)
 we obtain $K_2(\ab,x,y)=48(g^2-g+2)(c^2-4gk)x^2$ (respectively,
$K_2(\ab_1,x,y)=24g^2(c^2-8gk)x^2$). Calculations yield that for the system
(\ref{4s1a}) (respectively, for the system (\ref{5s4a})) we have
$\mu_0=\mu_1=0$ (respectively $\kappa=\kappa_1=0$). Hence we conclude that the
$GL$-comitant $K_2(a,x,y)$ is a  $CT$-comitant modulo
$\langle\eta,\mu_0,\mu_1,\kappa,\kappa_1\rangle$.

$d)$. We examine now  the comitant $K_3$ which is applied for systems
($\SSS_{I\!V}$) only in the cases when  $\Delta_S\ge3$, i.e.
$\mu_0=\mu_1=\mu_2=0$. It was shown before (see page \pageref{2s}) that  for
$\mu_0=\mu_1=\mu_2=0$ the systems $(\SSS_{I\!V})$ can be brought by an affine
transformation either to the systems (\ref{2s}) for $K\ne0$ or to the systems
(\ref{sys:K3}) for $K=0$. Calculations yield,  that for any system  in the
orbit under the translation group action of a system (\ref{2s}) (respectively,
of a system (\ref{sys:K3})) corresponding to a point $\ab\in \R^{12}$
(respectively, $\ab_1\in \R^{12}$)
 we obtain $K_3(\ab,x,y)=-12g^2lx^6$ (respectively,
$K_3(\ab_1,x,y)=3f(2c-f)x^6$). Hence in both cases the values of $K_3$ do not
depend of the vector defining the  translations. Therefore the $GL$-comitant
$K_3(a,x,y)$ is a  $CT$-comitant modulo $\langle M,\mu_0,\mu_1,\mu_2\rangle$.

{\it III) The cases 12--16}. Let $\tau\in T(2,\R)$ be the translation:
$x=\tilde x+\alpha$, $y=\tilde y+\beta$ and   consider a quadratic system
(\ref{2s1}) which corresponds to a point  $\ab\in \R^{12}$.  It is sufficient
to verify that  the following relations occur, where  $\xi=\tilde x\beta-
\tilde y \alpha$:
$$
\bal
& K_1(r_\tau\cdot a,\ \tilde x,\tilde y)= K_1(a,\tilde x,\tilde y) -
      \xi K(a,\tilde x,\tilde y);\\
& \mu_s(r_\tau\cdot \ab,\, \tilde x,\tilde y)= \mu_s(\ab,\,\tilde x,\tilde y)
+\sum_{k=0}^{s-1}{4-k\choose s-k}\xi^{s-k}\mu_k(\ab,\,\tilde x,\tilde y),\quad
s=1,2,3,4.
\eal
$$
 So, Lemma \ref{Table:Propreties} is proved and this  completes the proof of
the Theorem \ref{th_2}. \EProof

{\footnotesize

}
 \newpage
\appendix

\begin{center}{\bf APPENDIX}
\end{center}
 Let us consider the tensorial form of quadratic system:
$$
  \frac {dx^j}{dt}=a^j +a^j_{\alpha}x^{\alpha}
  + a^j_{\alpha\beta}x^{\alpha}x^{\beta}\quad (j,\alpha,\beta=1,2).
$$
The following invariants and comitants, defined by polynomials of $J_i,R_i$
which are tensorially defined $GL$-comitants, were used in \cite{Nik_Vlp} for
the classification in the neighbourhood of infinity of quadratic differential
systems:
$$
\begin{gathered}
2\tilde \mu= J_4,\ \ \tilde \sigma=  J_7,\ \ 2\tilde \theta = J_5,\ \  \tilde L= R_{12},\
\ 2\tilde M= 9R_3+ 6 R_6- 8 R_{11}^2,\ \  \tilde S_1= R_5,\\
\tilde S_2= 2J_1^2R_6+2J_1R_1^2 - 2J_2R_6 + J_2R_{11}^2+ 8J_3R_3 -8J_3R_6
-4 R_7 - R_8, \ \  \tilde H=R_{13},\\
\tilde S_3= R_{12}^2(7 J_2-6J_1^2 -8 J_3)- R_{12}(10 J_1R_5 +4 R_1R_{10}- 6
R_3R_9)
    + 4 R_3R_{10}^2 -4 R_5^2,\\
   \tilde S_4= 4J_3-J_2,\ \ \tilde V =  R_4^2 - R_2R_5,\ \  2\tilde A= 2 R_6 -3 R_3,\ \ 2\tilde \eta= J_4+ 20J_5- 8J_6\\
 2\tilde N= R_3, \ \ 2\tilde G= 2R_1^2 -2J_2R_3+ 4R_7+ R_8,\ \
2\tilde F=  J_2R_5 + 4 R_2R_3+ 4R_1R_4, \\
\end{gathered}
$$
where
$$
\begin{gathered}
J_{1} =
 a^\alpha_\alpha,\ \
J_{2}  =  a^\alpha_p a^\beta_q \ve_{\alpha \beta} \ve^{pq},\ \ J_{3} =
 a^\alpha a^\beta_{\alpha \beta},\ \
J_{4} = a^\alpha_{p r}a^\beta_{q k}a^\gamma_{s n}a^\delta_{l m}
        \ve_{\alpha \beta}\ve_{\gamma \delta}\ve^{pq}\ve^{rs}\ve^{kl}\ve^{mn},\\
\ J_{5} = a^\alpha_{\gamma }a^\beta_{\delta r}a^\gamma_{q k}a^\delta_{s l}
   \ve_{\alpha \beta}\ve^{pq}\ve^{rs}\ve^{kl},\ \
J_{6}= a^\alpha_{p r}a^\beta_{\alpha q}a^\gamma_{\delta s}
         a^\delta_{\beta \gamma}
\ve^{pq}\ve^{rs},\ \ J_{7}= a^\alpha_p a^\beta_{\gamma q}a^\gamma_{\alpha
\beta}\ve^{pq},\\
R_{1}= x^\alpha  a^\beta_q a^\gamma_{p \alpha} \ve_{\beta \gamma}\ve^{pq},\ \,
R_{2}= x^\alpha a^\beta a^\gamma_\alpha\ve_{\beta \gamma},\ \, R_{3}=
 x^\alpha x^\beta a^\gamma_{p \alpha}a^\delta_{q \beta} \ve_{\gamma \delta}\ve^{pq},\
 \, R_{4}= x^\alpha x^\beta a^\gamma
a^\delta_{\alpha \beta}\ve_{\gamma \delta},\\
R_{5} =
 x^\alpha x^\beta x^\gamma
 a^\delta_\alpha
a^\mu_{\beta \gamma}
        \ve_{\delta \mu},\ \
R_{6}= x^\alpha x^\beta a^\gamma_{\alpha \beta}a^\delta_{\gamma \delta},\ \
R_{7} =
 x^\alpha x^\beta
 a^\gamma
a^\delta_{\alpha p}a^\mu_{\beta s}
          a^\nu_{q r}
\ve_{\gamma \delta}\ve_{\mu \nu}\ve^{pq}\ve^{rs},\\
R_{8}=
 x^\alpha x^\beta
a^\gamma_\alpha a^\delta_\beta a^\mu_{p r}
          a^\nu_{q s}\ve_{\gamma \mu}\ve_{\delta \nu}\ve^{pq}\ve^{rs},\ \
R_{9}= x^\alpha
 a^\beta
\ve_{\beta \alpha},\ \ R_{10}= x^\alpha x^\beta
 a^\gamma_\alpha  \ve_{\gamma \beta},\\
R_{11}= x^\alpha a^\beta_{\alpha \beta},\ \ R_{12} =  x^\alpha x^\beta x^\gamma
a^\delta_{\alpha \beta}
         \ve_{\delta \gamma},\ \
R_{13}= x^\alpha
 a^\beta_p
a^\gamma_{\alpha r}a^\delta_{q k}a^\mu_{s l}
     \ve_{\beta \gamma}\ve_{\delta \mu}
\ve^{pq}\ve^{rs}\ve^{kl},\\
\end{gathered}
$$
 and
$$
\ve^{11}= \ve^{22}=\ve_{11}=\ve_{22}=0, \quad
\ve^{12}=\ve_{12}=-\ve^{21}=-\ve_{21}=1.
$$
\end{document}